\documentstyle[11pt,amsfonts]{article}
\title{Closed surfaces with bounds on their Willmore energy}
\author{Ernst Kuwert \& Reiner Sch\"atzle\footnote{Both authors
were supported by the Deutsche Forschungsgemeinschaft 
via DFG Forschergruppe 469, and by the Centro di Ricerca Matematica
Ennio De Giorgi during a visit in Pisa.}}

\newcommand {\rel} {{\mathbb R}}
\newcommand {\com} {{\mathbb C}}
\newcommand {\nat} {{\mathbb N}}

\newcommand {\ganz} {{\mathbb Z}}

\newtheorem{proposition}{Proposition}[section]
\newtheorem{theorem}{Theorem}[section]

\newtheorem{lemma}[proposition]{Lemma}

\newcommand{\bcite}[1] {\cite{#1}}

\newcommand {\pr} {\bf}

\newcommand{\proof} {   \begin{flushright}
                        ///
                        \end{flushright}
                }

\newcommand{\defin} { \hspace*{\fill} $\Box$ }

\newcommand{\mint}[1]{\mbox{$\displaystyle{
-\hspace{-1.05em}\int_{#1}} $}}

\newcommand {\D}{\displaystyle}

\def    \har    {{ {\cal H}^1 }}

\def    \mean   {{\vec H }}

\def    \Lt	{{ {\cal L}^2 }}
\def    \Lo	{{ {\cal L}^1 }}
\def	\d	{{ \ {\rm d} }}
\def	\ds	{{ {\rm d} }}

\def	\pro	{{ \pi }}
\def	\prog	{{ \Pi }}
\def    \ve     {{ \varepsilon }}

\def    \W      {{ {\mathcal W} }}
\def    \C      {{ {\mathcal C} }}

\def	\sur	{{ \Sigma }}

\def	\cor	{{ \cal P }}
\def	\grass	{{ \cal N }}
\def    \spt    {{\rm spt\,}}

\newcommand{\dd}{ \begin{displaymath} }
\newcommand{\df}{ \end{displaymath} }
\newcommand{\dcd}{ \begin{displaymath} \begin{array}{c}}
\newcommand{\dcf}{ \end{array} \end{displaymath} }
\newcommand{\ee}{ \begin{equation} }
\newcommand{\ef}{ \end{equation} }
\newcommand{\ad}{ \begin{array}{c} }
\newcommand{\af}{ \end{array} }

\begin{document}
\maketitle
\begin{abstract}\noindent{The Willmore energy of a closed
surface in $\rel^n$ is the integral of its squared mean 
curvature, and is invariant under M\"obius transformations 
of $\rel^n$. We show that any torus in $\rel^3$ with 
energy at most $8\pi-\delta$ has a representative under 
the M\"obius action for which the induced metric and a 
conformal metric of constant (zero) curvature are 
uniformly equivalent, with constants depending only on 
$\delta >0$. An analogous estimate is also obtained
for surfaces of fixed genus $p \geq 1$ in $\rel^3$ or $\rel^4$, 
assuming suitable energy bounds which are sharp for $n=3$. 
Moreover, the conformal type is controlled in terms of 
the energy bounds.
}
\end{abstract}

\setcounter{equation}{0}
\section{Introduction} \label{intro}

For an immersed surface $f : \sur \rightarrow \rel^n$ the Willmore functional
is defined as the integral
\begin{displaymath}
	\W(f) = \frac{1}{4} \int_\sur |\mean|^2 \d \mu_g,
\end{displaymath}
where $\mean$ is the mean curvature vector, 
$g = f^* g_{euc}\ $ is the pull-back metric
and $\mu_g$ is the induced area measure on $\sur$. 
The Gau{\ss} equation says that 
\begin{equation}
\label{intro.gauss1}
K = \frac{1}{2} (|\mean|^2 - |A|^2) = 
\frac{1}{4} |\mean|^2 - \frac{1}{2} |A^\circ|^2, 
\end{equation}
where $A_{ij} = A^\circ_{ij} + \frac{1}{2}\mean g_{ij}$ is the vector-valued 
second fundamental form and $K$ is the sectional curvature of $g$. In the case when
$\sur$ is an oriented closed surface of genus $p$, the Gau{\ss}-Bonnet theorem 
therefore implies the identities
\ee \label{intro.gauss}
	\W(f) = \frac{1}{4} \int_\sur
	|A|^2 \d \mu_g + 2 \pi (1 - p)
	= \frac{1}{2} \int_\Sigma
	|A^\circ|^2 \d \mu_g + 4 \pi (1 - p).
\ef
We denote by $\beta_p^n$ the infimum of the Willmore functional
among closed oriented surfaces $f:\sur \to \rel^n$ of genus $p$.
It is well-known that $\beta^n_0 = 4\pi$ with round spheres as
unique minimizers. For $p \geq 1$ we have $\beta^n_p > 4\pi$ 
by the analysis of L.\,Simon \cite{sim.will}. We put
\begin{equation}
\label{eqdouglasnumber}
\tilde{\beta}^n_p = \min \Big\{
4\pi + \sum_{i=1}^k (\beta^n_{p_i}-4\pi): 1 \leq p_i < p,\,\sum_{i=1}^k p_i =p\Big\},
\end{equation}
where $\tilde{\beta}^n_1 = \infty$, and define the constants
\begin{equation}
\label{eqomegaconstants}
\omega^n_p =
\left\{ \begin{array}{ll}
\min(8\pi,\tilde{\beta}^3_p) & \mbox{ for } n=3,\\
\min(8\pi,\tilde{\beta}^4_p,\beta^4_p + \frac{8\pi}{3}) &  \mbox{ for } n=4.
\end{array}
\right.\vspace{3mm}
\end{equation}
The main result of this paper is the following bilipschitz estimate. 
As the Willmore functional is invariant under the M\"obius group
of $\rel^n$, i.e. under dilations and inversions, the choice
of the M\"obius transformation in the statement is essential.

{\flushleft {\bf Theorem 4.1 }}{\em For $n =3,4$ and $p \geq 1$, let 
$\C(n,p,\delta)$ be the class of closed, oriented, genus $p$ surfaces
$f: \sur \rightarrow \rel^n$ satisfying $\W(f) \leq \omega^n_p - \delta$
for some $\delta > 0$. Then for any  $f \in \C(n,p,\delta)$ there is a
M\"obius transformation $\phi$ and a constant curvature metric $g_0$, 
such that the metric $g$ induced by $\phi \circ f$ satisfies}
$$
g = e^{2u}g_0 \quad \mbox{ where }\quad 
\max_{\sur} |u| \leq C(p,\delta) < \infty.
$$
We have $\beta^n_p < 8\pi$ as observed by Pinkall and independently Kusner, see for 
example \cite{kusner}, and $\beta^n_p < \tilde{\beta}^n_p$ from \cite{bauer.kuwert}.
Thus $\C(n,p,\delta)$ is nonempty at least for small $\delta >0$.
The stereographic projection of the Clifford torus into $\rel^3$ 
has energy $2\pi^2 < 8\pi = \omega^3_1$
and is conjectured to be the minimizer for $p=1$, compare \cite{schmidt1}.
We remark that we would have $\omega^n_p = 8\pi$ once we knew that
$\beta^n_{q} \geq 6\pi$ for $1 \leq q <p$ and $\beta^4_p \geq 16 \pi/3$
for $n=4$. It will be shown that our energy assumptions are sharp for $n=3$, 
that is the conclusion of the theorem fails if $\omega^3_p$ is replaced 
by any bigger constant. Combining the estimate in Theorem \ref{moeb.theo} 
with Mumford's compactness lemma we prove the following application.

{\flushleft {\bf Theorem 5.1 }}{\em For $n \in \{3,4\}$ and $p \geq 1$, 
the conformal structures induced by immersions $f$ in $\C(n,p,\delta)$ 
are contained in a compact subset $K = K(p,\delta)$ of the moduli space.}\\
\\
In particular as $\omega^4_1 \geq 20\pi/3$ we conclude 
$\W(f) > 2\pi^2$ for all tori $f:\sur \to \rel^4$ whose conformal structure is 
sufficiently degenerate.
A straightforward second application of Theorem 4.1 is a 
compactness theorem, which will be stated in our forthcoming paper 
\cite{kuw.schae.will4}. There the problem of minimizing the Willmore 
functional with prescribed conformal type is addressed.\\
\\
We now briefly summarize the contents of the paper. In Section 2 we review 
a version of the approximate graphical decomposition lemma on annuli, 
due to L. Simon \cite{sim.will}. In Section 3 we present the global estimate of the 
conformal factor, under certain technical assumptions. The choice of the 
M\"obius transformation is carried out in Section 4, and the proof of
Theorem 4.1 is then completed by verifying the assumptions from Section 3. 
In Section 5 we discuss the bound for the conformal type and the optimality 
of the constant $\omega^n_p$. Our results rely on estimates for surfaces of the 
type of the plane due to S. M\"uller and V. \v{S}verak \cite{muell.sver}.
The version needed is presented in the final section 6.\\
\\
{\bf Acknowledgement: }The main ideas of this paper were developed during a 
joint visit at the Centro di Ricerca Matematica Ennio De Giorgi, Pisa.
It is a pleasure to thank for the hospitality and the fruitful
scientific atmosphere.

\setcounter{equation}{0}

\section{Preliminaries} \label{pre}

Here we collect some results from the work of L.\,Simon \cite{sim.will},
starting with consequences of the monotonicity identity. 
For a proper immersion $f:\sur \to B_\varrho(0) \subseteq \rel^n$ of an 
open surface $\sur$ and any $\sigma \in (0,\varrho)$, we have by (1.3) 
in \cite{sim.will} the bound
\ee 
\label{pre.li-yau.mono}
\sigma^{-2} \mu(B_\sigma(0)) \leq 
C\,\Big(\varrho^{-2} \mu(B_\varrho(0)) + {\cal W}(f,B_\varrho(0))\Big),
\ef
where $\mu = f(\mu_g)$ is the pushforward area measure and
$$
{\cal W}(f,B_\varrho(0)) = \frac{1}{4}\,\int_{B_\varrho(0)} |\mean|^2 \d \mu. 
$$
We should really integrate over $f^{-1}(B_\varrho(0))$ with respect
to $\mu_g$, but the pullback is omitted for convenience;
in fact the notation can be justified by considering $\mu$ as a $2$-varifold 
with square integrable weak mean curvature as in the appendix of 
\cite{kuw.schae.will3}. If $\sur$ is compact without boundary we may let
$\varrho \nearrow \infty$ in (\ref{pre.li-yau.mono}) to get
\ee \label{pre.li-yau.dense}
\sigma^{-2} \mu(B_\sigma(0)) \leq C\, {\cal W}(f)
\quad \mbox{for all } \sigma > 0.
\ef
Moreover, the multiplicity of the immersion at $0$ is just the $2$-density
of $\mu$ and satisfies the Li-Yau inequality, see Theorem 6 in \cite{li.yau},
\ee \label{pre.li-yau}
\theta^2(\mu,0) \leq \frac{1}{4 \pi} {\cal W}(f).
\ef
We will need the following version of the approximate graphical decomposition 
lemma, see Lemma 2.1 and pp. 312--315 in \cite{sim.will}. 

\begin{lemma} \label{lemma.pre.graph} For any $\Lambda < \infty$
there exist $\ve_0 = \ve_0(n,\Lambda) > 0$ and $C = C(n,\Lambda) < \infty$ 
such that if $f:\sur \to B_\varrho(0) \subseteq \rel^n$ 
is a proper immersion satisfying
\begin{equation}
\label{pre.bounds}
\mu(B_\varrho(0) - B_{\varrho/2}(0)) \leq \Lambda \varrho^2
\quad \mbox{ and } \quad
\int_{B_\varrho(0) - B_{\varrho/2}(0)} |A|^2 \d \mu \leq \varepsilon^2
\quad \mbox{ for }\ve < \ve_0,
\end{equation}
then the following statements hold:
\begin{itemize}
\item[{\rm (a)}] Denote by $A^i$, $i=1,\ldots,m$, those components of
$f^{-1}\big(B_{7\varrho/8}(0)-B_{5\varrho/8}(0)\big)$ which extend to 
$\partial B_{9\varrho/16}(0)$. There exist compact subdiscs 
$P_1,\ldots,P_N \subseteq \sur$ with
$$
\sum_{j=1}^N {\rm diam\,}f(P_j) < C \ve^{1/2} \varrho,
$$
such that on each $A^i - \bigcup_{j=1}^N P_j$ the immersion is 
a $k_i$-valued graph for $k_i \in \nat$, intersected with 
$B_{7\varrho/8}(0)-B_{5\varrho/8}(0)$, over some affine $2$-plane.
Furthermore 
\begin{equation}
\label{pre.graph.mult}
M: = \sum_{i=1}^m k_i \leq C.
\end{equation}
\item[{\rm (b)}] There is a set $S \subseteq (5\varrho/8,7\varrho/8)$ of measure
$\Lo(S) > 3\varrho/16$, such that for $\sigma \in S$ the immersion is transversal
to $\partial B_\sigma(0)$, 
and for each $\Gamma^i_\sigma: = A^i \cap f^{-1}\big(\partial B_\sigma(0)\big)$ 
we have 
\begin{equation}
\label{eqgeodesiccurvaturebound}
\Big|\int_{\Gamma^i_\sigma} \kappa_g\,\ds s - 2\pi k_i\Big| \leq C \ve^\alpha
\quad \mbox{ where }\alpha = \alpha(n) > 0.
\end{equation}
Furthermore, the restriction of $f$ to $A^i_\sigma:= A^i \cap f^{-1}(B_\sigma(0))$ 
has a $C^{1,1}$ extension $\tilde{f}:\tilde{A}^i_\sigma \to \rel^n$, where
$\tilde{A}^i_\sigma$ is obtained by attaching a punctured disc $E^i_\sigma$ 
to $A^i_\sigma$ along $\Gamma^i_\sigma$, such that
$\tilde{f}$ is a flat $k_i$-fold covering of an affine $2$-plane $L_i$ outside 
$B_{2\sigma}(0)$ and has curvature bounded by
\begin{equation} \label{eqcurvaturebound}
\int_{E^i_\sigma} |\tilde{A}|^2 \d\tilde{\mu} \leq C \ve^2.
\end{equation}
\item[{\rm (c)}] If we assume in addition to {\rm (\ref{pre.bounds})} that 
\begin{equation}
\label{pre.heightestimate}
\int_{B_\varrho(0) - B_{\varrho/2}(0)}
\frac{|x^{\perp}|^2}{|x|^4} \d \mu(x) < \varepsilon^2,
\end{equation}
where $\perp$ denotes the projection in the normal direction along the immersion,
then $f^{-1}\big(B_{7\varrho/8}(0)-B_{5\varrho/8}(0)\big) = \bigcup_{i=1}^m A^i$
and we have the estimate
\begin{equation}
\label{pre.graph.area}
\mu(B_{7\varrho/8}(0)-B_{5\varrho/8}(0)) \geq (1-C\varepsilon^{2\alpha})
M\pi\,\left((7\varrho/8)^2 - (5\varrho/8)^2 \right).
\end{equation}
\end{itemize}
\end{lemma}
If the assumptions of Lemma \ref{lemma.pre.graph} hold with $\varrho/2$ replaced by 
some $r \in (0,\varrho/2]$, that is $\mu\big(B_\varrho(0)-B_r(0)\big) \leq \Lambda \varrho^2$ 
and 
\begin{equation}
\label{pre.graph.densityass}
\int_{B_\varrho(0) - B_r(0)} |A|^2 \d \mu,\,
\int_{B_\varrho(0) - B_r(0)} \frac{|x^\perp|^2}{|x|^4}\d \mu(x) < \varepsilon^2,
\end{equation}
then by inequality (\ref{pre.li-yau.mono}) the assumptions of Lemma
\ref{lemma.pre.graph} are satisfied with $\varrho$ replaced by any 
$\sigma \in [2r,\varrho]$. The resulting graphical decompositions
have the same multiplicity $M$ by continuity. Choosing $\sigma_\nu = (5/7)^\nu \varrho$
and summing over the inequalities (\ref{pre.graph.area}) we find
\begin{equation}
\label{pre.graph.density}
\mu(B_{7\varrho/8}(0)-B_{5r/4}(0)) \geq (1-C\varepsilon^{2\alpha})
M\pi\,\left((7\varrho/8)^2 - (5r/4)^2 \right).
\end{equation}
The results in \cite{sim.will} are stated only for
embedded surfaces, however they extend to immersions simply by considering
a pertubation $f_\lambda = (f,\lambda f_0):\sur \to \rel^n \times \rel^3$,
where $f_0:\sur \to \rel^3$ is some differentiable embedding. The
$f_\lambda$ satisfy the assumptions of Lemma \ref{lemma.pre.graph}
for a slightly bigger constant $\Lambda$, hence they admit a graphical
decomposition as stated over some $2$-planes in $\rel^n \times \rel^3$,
which are almost horizontal for $\lambda$ sufficiently small.
By slightly tilting the planes one obtains the desired almost
graphical decomposition for the given immersion $f$, with power $\alpha = 1/(4n+6)$
instead of $1/(4n-6)$ which is the constant in Lemma 2.1 of \cite{sim.will}.

\setcounter{equation}{0}

\section{Oscillation estimates} \label{osc}

In this section we present the main PDE argument for the estimate of the 
conformal factor. 

\begin{theorem} \label{osc.theo}

Let $f: \sur \rightarrow \rel^n$, $n = 3,4$, be an immersion
of a closed surface $\sur$ of genus $p \geq 1$
with ${\cal W}(f) \leq \Lambda$. Assume that 
$f(\sur) \subseteq \bigcup_{k=1}^K B_{\varrho_k/2}(x_k)$ with
$\varrho_l/\varrho_k \leq \Lambda$, such that for all $k=1,\ldots,K$ 
and some $\delta > 0$ the following conditions hold:
\begin{equation} \label{osc.theo.curv3}
	\int_{B_{\varrho_k}(x_k)}
	|K| \d \mu < 8 \pi - \delta
	\quad \mbox{ for } n = 3,
\end{equation}
\begin{equation} \label{osc.theo.curv4}
\left.
\begin{array}{ccl}
	\displaystyle{\int_{B_{\varrho_k}(x_k)} |K| \d \mu}
	+ \displaystyle{\frac{1}{2} \int_{B_{\varrho_k}(x_k)} |A^\circ|^2 \d \mu}
	& < & 8 \pi - \delta\\
	\displaystyle{\int_{B_{\varrho_k}(x_k)} |A^\circ|^2 \d \mu}
	& \leq  & 8 \pi - C_0 \varepsilon^2
\end{array}
\right\}
	\quad \mbox{for } n = 4,
\end{equation}
\begin{equation} \label{osc.theo.curv-ann}
	\int_{B_{\varrho_k}(x_k) - B_{\varrho_k/2}(x_k)}
	|A|^2 \d \mu < \varepsilon^2.
\end{equation}
Denoting by $D^{k,\alpha}_\sigma$, $1 \leq \alpha \leq m_k$, the 
components of $f^{-1}(B_{\sigma}(x_k))$ which meet 
$\partial B_{9\varrho_k/16}(x_k)$, we further assume for all
$\sigma \in [5\varrho_k/8,7\varrho_k/8]$ up to a set of 
measure at most $\varrho_k/16$ that 
\begin{equation}
\label{osc.theo.topo}
\int_{D^{k,\alpha}_\sigma} K \d\mu_g > -2\pi + \delta \quad
\mbox{ for all }\alpha = 1,\ldots,m_k.
\end{equation}
Then for $\varepsilon \leq \varepsilon(\Lambda,\delta)$ and $C_0 \geq C_0(\Lambda)$, 
there is a constant curvature metric $g_0 = e^{-2u}g$ such that
\begin{displaymath}
        \max_{\sur} |u| \leq C(\Lambda,K,p,\delta).
\end{displaymath}
\end{theorem}
{\pr Proof:}
By rescaling we may assume $\mu_g(\sur) = 1$. We take $g_0 = e^{-2u}g$ 
as the unique conformal, constant curvature metric also with 
$\mu_{g_0}(\sur) = 1$, which means
\begin{equation}
\label{eqconformalfactor}
-\Delta_g u + K_{g_0}\,e^{-2u} = K_g \quad \mbox{ where }\quad
K_{g_0} = \frac{2 \pi \chi(\Sigma)}{\mu_{g_0}(\Sigma)} = 4\pi(1-p).
\end{equation}
Clearly the condition $\mu_{g_0}(\sur) = \mu_g(\sur)$ implies 
\begin{equation}
\label{eqconffactorzero}
u(p)=0 \quad \mbox{ for some } p \in \sur,
\end{equation}
and hence it suffices to prove the estimate
\ee \label{osc.theo.aux}
   {\rm osc}_\sur u \leq C(\Lambda,K,p,\delta).
\ef
The bound ${\cal W}(f) \leq \Lambda$ and the identity (\ref{intro.gauss}) imply
\ee \label{osc.theo.funda}
	\int_\sur |A|^2 \d \mu_g
	\leq C(\Lambda,p),
\ef
and the Li-Yau inequality (\ref{pre.li-yau.dense}) yields
\ee \label{osc.theo.dense-ratio}
	\varrho^{-2} \mu(B_\varrho(x)) \leq C(\Lambda)
	\quad \mbox{for all } B_\varrho(x) \subseteq \rel^n.
\ef
The set of all $\sigma \in [5\varrho_k/8,7\varrho_k/8]$ 
satisfying both (\ref{osc.theo.topo}) and the inequality
\begin{equation}
\label{pre.goodslice.ass}
\int_{\partial B_\sigma(x_k)} |A|^2\,\d s:=
\int_{\partial \big[f^{-1}(B_\sigma(x_k))\big]} |A|^2\,\d s_g
\leq  16 \varepsilon^2/\varrho_k
\end{equation}
has measure at least $\varrho_k/8$. Thus we can choose
$\sigma_k,\sigma'_k \in [5\varrho_k/8,7\varrho_k/8]$ satisfying 
(\ref{osc.theo.topo}), (\ref{pre.goodslice.ass}) and the
conclusions of Lemma \ref{lemma.pre.graph}(b), such that
$\sigma_k - \sigma'_k > \varrho_k/16$. Since 
$f(\sur) \subseteq \bigcup_{k=1}^K B_{\varrho_k/2}(x_k)$ we have
\ee \label{osc.theo.cover}
        \sur = \bigcup_{k,\alpha} D^{k,\alpha}_{\sigma'_k}.
\ef
From the Gau{\ss}-Bonnet theorem and (\ref{osc.theo.topo}),
we obtain for each component
$$
\int_{\partial D^{k,\alpha}_{\sigma_k}} \kappa_g \,\ds s_g 
= 2\pi \chi(D^{k,\alpha}_{\sigma_k}) - \int_{D^{k,\alpha}_{\sigma_k}} K \d\mu_g 
< 2\pi\left(\chi(D^{k,\alpha}_{\sigma_k}) + 1\right)-\delta.
$$
We conclude that each $D^{k,\alpha}_{\sigma_k}$ is a disc, and that 
the multiplicity of its boundary entering in (\ref{eqgeodesiccurvaturebound})
equals one, which means that all the graphs in  Lemma \ref{lemma.pre.graph}(a) 
are singlevalued. Again by Lemma \ref{lemma.pre.graph}(b), we extend
$f|_{D^{k,\alpha}_{\sigma_k}}$ to an immersion 
$f_{k,\alpha}:\sur_{k,\alpha} \rightarrow \rel^n$ such that 
\ee \label{osc.theo.ext-funda}
	\int_{\sur_{k,\alpha} - D^{k,\alpha}_{\sigma_k}}
	|A_{f_{k,\alpha}}|^2 \d \mu_{f_{k,\alpha}}
	\leq C \varepsilon^2.
\ef
Here $D^{k,\alpha}_{\sigma_k} \subseteq \sur_{k,\alpha} \cong \rel^2$
and $f_{k,\alpha}$ is the standard embedding of a single plane 
outside $B_{2 \sigma_k}(x_k) \subseteq \rel^n$, in particular 
$f_{k,\alpha}$ is complete. Now for $g_{k,\alpha} := f_{k,\alpha}^* g_{euc}$
the Gau{\ss}-Bonnet theorem implies
\dd
	\int_{\sur_{k,\alpha}} K_{g_{k,\alpha}}
	\d \mu_{g_{k,\alpha}} = 0.
\df
By the uniformization theorem, we may assume that
the diffeomorphism $\sur_{k,\alpha} \cong \rel^2$ is conformal, and write
$g_{k,\alpha} = e^{2 u_{k,\alpha}} g_{euc}$ on $\sur_{k,\alpha} \cong \rel^2$.
From (\ref{osc.theo.curv3}), (\ref{osc.theo.curv4})
and (\ref{osc.theo.ext-funda}), we get
\dd
	\int_{\sur_{k,\alpha}}
	|K_{g_{k,\alpha}}| \d \mu_{g_{k,\alpha}}
	\leq 8 \pi - \delta + C \varepsilon^2
	\quad \mbox{for } n = 3,
\df
\dd
\left.
\ad\D
	\int_{\sur_{k,\alpha}}
	|K_{g_{k,\alpha}}| \d \mu_{g_{k,\alpha}}
	+ \frac{1}{2} \int_{\sur_{k,\alpha}}
	|A^\circ_{f_{k,\alpha}}|^2 \d \mu_{g_{k,\alpha}}
	\leq 8 \pi - \delta + C \varepsilon^2 \\
	\int_{\sur_{k,\alpha}}
	|A^\circ_{f_{k,\alpha}}|^2 \d \mu_{g_{k,\alpha}}
	\leq 8 \pi - C_0 \varepsilon^2 + C \varepsilon^2
\af
\right\}
	\quad \mbox{for } n = 4.
\df
Choosing $C \varepsilon^2 < \delta/2$ and $C_0 > C$,
this verifies the assumptions of Theorem \ref{conf.flat},
except that the parameter $\delta$ is replaced by $\delta/2$. Thus we have
\dd
	- \Delta_{g_{k,\alpha}} u_{k,\alpha} = K_{g_{k,\alpha}}
	\quad \mbox{in } \sur_{k,\alpha},
\df
where $u_{k,\alpha}$ satisfies the estimates, possibly after adding a suitable constant, 
\begin{eqnarray*}
\|u_{k,\alpha}\|_{L^\infty(\sur_{k,\alpha})},
\|D u_{k,\alpha}\|_{L^2(\sur_{k,\alpha})},
\|D^2 u_{k,\alpha}\|_{L^1(\sur_{k,\alpha})} & \leq &
C(\delta)\, \int_{\sur_{k,\alpha}}
	|A_{f_{k,\alpha}}|^2 \d \mu_{g_{k,\alpha}}\\
& \leq & C(\Lambda,p,\delta).
\end{eqnarray*}
Here the $L^1$ and $L^2$ norms on the left are with respect to the Euclidean metric
on $\sur_{k,\alpha} \cong \rel^2$, and we use (\ref{osc.theo.funda}) and 
(\ref{osc.theo.ext-funda}) for the second inequality.
As $f_{k,\alpha}$ and $f$ coincide on $D^{k,\alpha}_{\sigma_k}$,
we have $g_{k,\alpha} = g$ on $D^{k,\alpha}_{\sigma_k}$,
hence 
\ee \label{osc.theo.solu}
	- \Delta_g u_{k,\alpha} = K_g
	\quad \mbox{in } D^{k,\alpha}_{\sigma_k}.
\ef
and by conformal invariance of the Dirichlet integral
\ee \label{osc.theo.solu.esti}
\|u_{k,\alpha}\|_{L^\infty(D^{k,\alpha}_{\sigma_k})},
\int_{D^{k,\alpha}_{\sigma_k}} |D u_{k,\alpha}|_g^2 \d \mu_g
\leq C(\Lambda,p,\delta).
\ef
Combining with (\ref{eqconformalfactor})
and as $g = g_{k,\alpha} = e^{2 u_{k,\alpha}} g_{euc}$
on $D^{k,\alpha}_{\sigma_k} \subseteq \sur_{k,\alpha}
\cong \rel^2$, we get
\dd
	- \Delta (u - u_{k,\alpha})
	= - e^{2 u_{k,\alpha}} \Delta_g (u - u_{k,\alpha})
	= - K_{g_0} e^{-2(u - u_{k,\alpha})}
	\quad \mbox{ in } D^{k,\alpha}_{\sigma_k},
\df
hence using $0 \leq - K_{g_0} = 4 \pi (p-1)$ we conclude
\ee \label{osc.theo.equ-local}
\left.
\ad
	- \Delta(u - u_{k,\alpha}) \geq 0, \\
	- \Delta\,(u - u_{k,\alpha})_+ \leq C(p-1),	
\af
\right\}
	\quad \mbox{in } D^{k,\alpha}_{\sigma_k}.
\ef
Next we choose extrinsic cut-off functions
$\gamma_k \in C^2_0(B_{5 \varrho_k/8}(x_k))$ with $0 \leq \gamma_k \leq 1$,
$\gamma_k = 1$ on $B_{\varrho_k/2}(x_k)$ and 
$|D^j \gamma_k| \leq C \varrho_k^{-j}$ for $j = 1,2$; we then put
$\tilde \eta_{k,\alpha} := \gamma_k \circ f$ on $D^{k,\alpha}_{\sigma_k}$
and $\tilde \eta_{k,\alpha} = 0$ on $\sur - D^{k,\alpha}_{\sigma_k}$.
Then $\tilde \eta := \sum_{k,\alpha} \tilde \eta_{k,\alpha} \geq 1$
on $\sur$, as $f(\sur)$ is covered by the $B_{\varrho_k/2}(x_k)$ for 
$k=1,\ldots,K$. We put $\eta_{k,\alpha} =  \tilde \eta_{k,\alpha}/ \tilde \eta$ and get
\dcd
	\spt \eta_{k,\alpha} \subseteq D^{k,\alpha}_{\sigma_k}, \\
	\sum_{k,\alpha} \eta_{k,\alpha} = 1, \\
	|D \eta_{k,\alpha}|_g \leq C(K) \varrho_k^{-1}, \\
	|D^2 \eta_{k,\alpha}|_g \leq C(\Lambda,K)
	(\varrho_k^{-2} + \varrho_k^{-1} |A|).
\dcf
Putting $\bar u := \sum_{k,\alpha} \eta_{k,\alpha} u_{k,\alpha}$,
we calculate from (\ref{eqconformalfactor}) and (\ref{osc.theo.solu})
\ee \label{osc.theo.equ-global}
	- \Delta_g (u-\bar u)
	= - K_{g_0} e^{-2u} + 2 \sum_{k,\alpha}
	D \eta_{k,\alpha} D u_{k,\alpha}
	+ \sum_{k,\alpha} \Delta_g \eta_{k,\alpha}
	\ u_{k,\alpha} =: h
\ef
and estimate by (\ref{osc.theo.funda}), (\ref{osc.theo.dense-ratio}),
(\ref{osc.theo.solu.esti}),
recalling $K_{g_0} = 4 \pi (1-p) \leq 0$ and $g_0 = e^{-2u} g$,
\begin{eqnarray*}
\int_\sur |h| \d \mu_g & \leq &  \int_{\sur} (- K_{g_0}) e^{-2u} \d \mu_g\\
&   & +\, C(\Lambda,K)\, \sum_{k,\alpha} \left(\varrho_k^{-2} \mu_g(D^{k,\alpha}_{\sigma_k})\right)^{1/2}
       \Big( \int_{D^{k,\alpha}_{\sigma_k}}  |D u_{k,\alpha}|_g^2 \d \mu_g \Big)^{1/2}\\
&   & +\, \|u_{k,\alpha}\|_{L^\infty(D^{k,\alpha}_{\sigma_k})}
       \int_{D^{k,\alpha}_{\sigma_k}} \big(\varrho_k^{-2} + \varrho_k^{-1} |A|\big) \d \mu_g\\
& \leq &  C(\Lambda,K,p,\delta).
\end{eqnarray*}
Furthermore 
\dd
	\|\bar u \|_{L^\infty(\sur)},
	\int_\sur |D \bar u|_g^2 \d \mu_g
	\leq C(\Lambda,K,p,\delta).	
\df
Multiplying (\ref{osc.theo.equ-global}) by $u - \bar u - \lambda$ where
$\lambda \in \rel$ is arbitrary, we obtain
\begin{eqnarray*}
\int_\sur |D(u - \bar u)|_g^2 \d \mu_g & \leq &
\int_\sur |h|\ |u - \bar u - \lambda| \d \mu_g\\
& \leq & C(\Lambda,K,p,\delta)\, \|u - \bar u - \lambda\|_{L^\infty(\sur)}\\
& \leq &  C(\Lambda,K,p,\delta)\, \big(1 + \|u - \lambda \|_{L^\infty(\sur)}),
\end{eqnarray*}
hence
\begin{eqnarray}
\nonumber 
\int_\sur |D u|_g^2 \d \mu_g & \leq &  	
2 \int_\sur |D(u - \bar u)|_g^2 \d \mu_g + 2 \int_\sur |D \bar u|_g^2 \d \mu_g\\
\label{osc.theo.grad}
&\leq &  C(\Lambda,K,p,\delta)\, (1 + {\rm osc}_\sur u).
\end{eqnarray}
Recalling the choice of $\sigma_k,\sigma_k'$, we note 
$B^g_{\varrho_k/16}(z) \subseteq D^{k,\alpha}_{\sigma_k}$ for $z \in D^{k,\alpha}_{\sigma'_k}$,
where $B^g_{\varrho}(z)$ is the geodesic ball with respect to $g$. Writing $B^2_\varrho(z)$ 
for the Euclidean coordinate disc using $\sur_{k,\alpha} \cong \rel^2$, we see from
(\ref{osc.theo.solu.esti}) that 
$B^2_{2 c_0 \varrho_k}(z) \subseteq D^{k,\alpha}_{\sigma_k}$
for $c_0 = c_0(\Lambda,p,\delta) > 0$ small enough. 
Now by (\ref{osc.theo.cover}) any $z \in \sur$ belongs to some $D^{l,\beta}_{\sigma'_l}$,
hence by (\ref{osc.theo.solu.esti}) and $\mu_g(\sur) =1$ 
\dd
	\pi (c_0 \varrho_l)^2 = \Lt(B_{c_0 \varrho_l}(z))
	\leq C(\Lambda,p,\delta)\, \mu_g(D^{l,\beta}_{\sigma_l})
	\leq C(\Lambda,p,\delta),
\df
hence $\varrho_k \leq C(\Lambda,p,\delta)$ for all $k$ since
$\varrho_k/\varrho_l \leq \Lambda$ by assumption.
Further by (\ref{osc.theo.dense-ratio}) 
\dd
	1 = \mu_g(\sur)
	\leq \sum_{k=1}^K \mu(B_{\varrho_k/2}(x_k))
	\leq C(\Lambda) K \max_{1 \leq k \leq K} \varrho_k^2.
\df
Using again $\varrho_l/\varrho_k \leq \Lambda$ we see that
\ee \label{osc.theo.radius}
	c_0(\Lambda,K) \leq \varrho_k \leq C(\Lambda,p,\delta).
\ef
Next, (\ref{osc.theo.grad}) and the Poincar\'e inequality show that, 
for appropriate $\lambda_{k,\alpha,z} \in \rel$,
$$
(c_0 \varrho_k)^{-1} \|u - \lambda_{k,\alpha,z}\|_{L^2(B^2_{c_0 \varrho_k}(z))} \leq 
C\, \|Du\|_{L^2(B^2_{c_0 \varrho_k}(z))} \leq
C\,\big(1 +  \sqrt{{\rm osc}_\sur u}\big). 
$$
Select a maximal subset $\{ z_i \}_{i \in I} \subseteq D^{k,\alpha}_{\sigma'_k}$
with $B^2_{c_0 \varrho_k/4}(z_i)$ pairwise disjoint, whence
the $B^2_{c_0 \varrho_k/2}(z_i)$, $i \in I$, cover $D^{k,l}_{\sigma'_k}$.
As the $D^{k,\alpha}_{\sigma_k} \supseteq B^2_{c_0 \varrho_k/4}(z_i)$
are pairwise disjoint, we estimate the cardinality of $I$ by
\dd
	{\rm card}(I)\ \pi (c_0 \varrho_k/4)^2
	\leq \Lt(D^{k,\alpha}_{\sigma_k})
	\leq C(\Lambda,p,\delta) \mu_g(D^{k,\alpha}_{\sigma_k}) \leq
\df
\dd
	\leq C(\Lambda,p,\delta) \mu(B_{\varrho_k}(x_k))
	\leq C(\Lambda,p,\delta) \varrho_k^2,
\df
as $g = e^{2 u_{k,\alpha}} g_{euc}$, using (\ref{osc.theo.solu.esti}) 
and (\ref{osc.theo.dense-ratio}), hence
\dd
	{\rm card}(I) \leq C(\Lambda,p,\delta).
\df
If $B^2_{c_0 \varrho_k/2}(z_i) \cap B^2_{c_0 \varrho_k/2}(z_j) \neq \emptyset$,
then $\Lt(B^2_{c_0 \varrho_k}(z_i) \cap B^2_{c_0 \varrho_k}(z_j)) 
\geq \pi (c_0 \varrho_k/2)^2$ and
\begin{eqnarray*}
|\lambda_{k,\alpha,z_i} - \lambda_{k,\alpha,z_j}| & \leq &
C\,(c_0 \varrho_k)^{-1}\left(\|u - \lambda_{k,\alpha,z_i}\|_{L^2(B^2_{c_0 \varrho_k}(z_i))} +
                       \|u-\lambda_{k,\alpha,z_j}\|_{L^2(B^2_{c_0 \varrho_k}(z_j))}\right)\\
&\leq & C(\Lambda,K,p,\delta)\,(1 + \sqrt{{\rm osc}_\sur u}).
\end{eqnarray*}
As $D^{k,\alpha}_{\sigma'_k}$ is connected and covered by the $B^2_{c_0 \varrho_k/2}(z_i)$,
we find for $i,j \in I$ a chain $B^2_{c_0 \varrho_k/2}(z_{i_\nu})$, $\nu=1,\ldots,N$,
with $N \leq {\rm card}(I)$ and such that neighboring discs intersect. Thus
$$
|\lambda_{k,\alpha,z_i} - \lambda_{k,\alpha,z_j}|
\leq C(\Lambda,K,p,\delta)\, (1 + \sqrt{{\rm osc}_\sur u}) \quad \forall i,j \in I.
$$
Therefore 
there exists a $\lambda_{k,\alpha} \in \rel$ such that
\dd
\varrho_k^{-1} \|u - \lambda_{k,\alpha}\|_{L^2(D^{k,\alpha}_{\sigma'_k};g)}
\leq C(\Lambda,K,p,\delta)\,(1 + \sqrt{{\rm osc}_\sur u}).
\df
The sets $B^{k,\alpha}_{\varrho_k/2}: = D^{k,\alpha}_{\sigma'_k} \cap f^{-1}(B_{\varrho_k/2}(x_k))$
form an open cover of $\sur$.  
Moreover if $z \in B^{k,\alpha}_{\varrho_k/2} \cap B^{l,\beta}_{\varrho_l/2}$
where $\varrho_k \leq \varrho_l$, then we obtain as above
$B^g_{\varrho_k/8}(z) \subseteq D^{k,\alpha}_{\sigma'_k} \cap D^{l,\beta}_{\sigma'_l}$,
using $\sigma'_k \geq 5 \varrho_k / 8$, $\sigma'_l \geq 5 \varrho_l / 8\ $, and
\dd
	\mu_g(D^{k,\alpha}_{\sigma'_k}
	\cap D^{l,\beta}_{\sigma'_l})
	\geq \mu_g(B^g_{\varrho_k / 8}(z))
	\geq c_0(\Lambda,p,\delta)
	\Lt(B^2_{c_0 \varrho_k}(z))
	\geq c_0(\Lambda,p,\delta) \varrho_k^2.
\df
This yields
\begin{eqnarray*}
|\lambda_{k,\alpha} - \lambda_{l,\beta}| & \leq &
(c_0 \varrho_k)^{-1}
\|\lambda_{k,\alpha} - \lambda_{l,\beta}\|_{L^2(D^{k,\alpha}_{\sigma'_k} \cap D^{l,\beta}_{\sigma'_l};g)}\\
&\leq & (c_0 \varrho_k)^{-1}\,\left(
	\|u - \lambda_{k,\alpha}\|_{L^2(D^{k,\alpha}_{\sigma'_k};g)}
      + \|u - \lambda_{l,\beta}\|_{L^2(D^{l,\beta}_{\sigma'_l};g)}\right)\\
&\leq & C(\Lambda,K,p,\delta)\, (1 + \sqrt{{\rm osc}_\sur u}),
\end{eqnarray*}
as $\varrho_l / \varrho_k \leq \Lambda$ by assumption. 
Again by connectedness of $\sur$ there is a $\lambda \in \rel$ 
\ee \label{osc.theo.mean}
\|u - \lambda \|_{L^2(\sur;g)} 
\leq C(\Lambda,K,p,\delta)\, (1 + \sqrt{{\rm osc}_\sur u}) 
\max_{1 \leq k \leq K} \varrho_k.
\ef
Next choose $z_0 \in \sur$ with $u(z_0) = \min_\sur u$.
Then $z_0 \in B^{k,\alpha}_{\varrho_k/2}$ for some $k,\alpha$. By
(\ref{osc.theo.solu.esti}) and (\ref{osc.theo.equ-local}) we have,
as $B^2_{2 c_0 \varrho_k}(z_0) \subseteq D^{k,\alpha}_{\sigma_k}$, the estimate
$u - u_{k,\alpha} \geq \min_\sur u - C(\Lambda,p,\delta)=: \bar \lambda$,
and conclude from the weak Harnack inequality, see \bcite{gil.tru} Theorem 8.18,
\dd
(c_0 \varrho_k)^{-1} \|u - u_{k,\alpha}- \bar \lambda \|_{L^2(B^2_{c_0 \varrho_k}(z_0))}
\leq C \inf_{B^2_{c_0 \varrho_k}(z_0)}  (u - u_{k,\alpha} - \bar \lambda).
\df
Hence from $u(z_0) = \min_\sur u$ we see that
\dd
(c_0 \varrho_k)^{-1} \|u - \min_\sur u\|_{L^2(B^2_{c_0 \varrho_k}(z_0))} \leq C(\Lambda,p,\delta).
\df
Then
\dd
|\min_\sur u - \lambda| \leq 
C(\Lambda,K,p,\delta)\, (1 + \sqrt{{\rm osc}_\sur u})  \max_{1 \leq k \leq K} \varrho_k
\df
by (\ref{osc.theo.mean}), and we conclude
\ee \label{osc.theo.mean-min}
\|u - \min_\sur u \|_{L^2(\sur,g)}  \leq 
C(\Lambda,K,p,\delta)\, (1 + \sqrt{{\rm osc}_\sur u})  \max_{1 \leq k \leq K} \varrho_k.
\ef
Now $\min_{\sur} u \leq 0$ by (\ref{eqconffactorzero}).
Employing the mean value inequality, see  \bcite{gil.tru} Theorem 2.1, we 
obtain from (\ref{osc.theo.equ-local}) for $z \in D^{k,\alpha}_{\sigma'_k}$
\dd
\|(u - u_{k,\alpha})_{+}\|_{L^\infty(B^2_{c_0 \varrho_k}(z))} \leq
C (c_0 \varrho_k)^{-1} \|(u - u_{k,\alpha})_+\|_{L^2(B^2_{2 c_0 \varrho_k}(z))}
+ C (c_0 \varrho_k)^2 (p-1).
\df
Combining (\ref{osc.theo.solu.esti}), (\ref{osc.theo.radius}),
(\ref{osc.theo.mean-min}) and $\min_{\sur} u \leq 0$ implies
\begin{eqnarray*}
\max_\sur u & \leq &
C(\Lambda,p,\delta) \varrho_k^{-1}
\|u - \min_\sur u \|_{L^2(\sur;g)} + C(\Lambda,p,\delta)+  C (c_0 \varrho_k)^2 (p-1)\\
&\leq & C(\Lambda,K,p,\delta)\,(1 + \sqrt{{\rm osc}_\sur u})\\
&\leq & C(\Lambda,K,p,\delta)  + \frac{1}{2} \max_\sur u - \frac{1}{2} \min_\sur u,
\end{eqnarray*}
hence $\max_\sur u \leq C(\Lambda,K,p,\delta) + |\min_\sur u|$, and
\ee \label{osc.theo.osc-min}
{\rm osc}_\sur u \leq C(\Lambda,K,p,\delta) + 2 |\min_\sur u|.
\ef
Next we define $A = \{x \in \sur: u(x) \leq \min_\sur u/2\}$.
As $u - \min_\sur u \geq |\min_\sur u|/2$ on $ \sur - A$,
we get from (\ref{osc.theo.radius}), (\ref{osc.theo.mean-min})
and (\ref{osc.theo.osc-min})
\begin{eqnarray*}
\frac{1}{2} |\min_\sur u|\, \mu_g(\sur - A) & \leq &
\int_{\sur} (u - \min_\sur u) \d \mu_g\\
&\leq & C(\Lambda,K,p,\delta)\, (1 + \sqrt{{\rm osc}_\sur u})\\
&\leq & C(\Lambda,K,p,\delta)\, (1 + \sqrt{|\min_\sur u|}).
\end{eqnarray*}
Thus for $|\min_\sur u| \gg C(\Lambda,K,p,\delta)$ we estimate
\dd
\mu_g(\sur - A)
\leq \frac{C(\Lambda,K,p,\delta)
(1 + \sqrt{|\min_\sur u|})}
{|\min_\sur u|} \leq \frac{1}{2}.
\df
As both $g$ and $g_0 = e^{-2u} g$ have unit area, this yields
$\mu_g(A) \geq 1/2$ and
\dd
1 \geq \int_A e^{-2u} \d \mu_g
\geq \mu_g(A)\, e^{-\min_\sur u}
\geq \frac{1}{2}\,e^{|\min_\sur u|}.
\df
We conclude that $|\min_\sur u| \leq C(\Lambda,K,p,\delta)$,
and hence (\ref{osc.theo.aux}) follows from (\ref{osc.theo.osc-min}),
and the theorem is proved. \proof
Inspecting the proof, we see that instead of (\ref{osc.theo.topo})
we could require directly that each component $D^{k,\alpha}_\sigma$ 
is a disc and $f|_{\partial D^{k,\alpha}_\sigma}$ is a single, nearly 
flat circle, for all $\sigma \in ]5\varrho_k/8,7\varrho_k/8[$ up to a 
set of measure $\varrho_k/16$. We also remark that the
assumptions (\ref{osc.theo.curv3})-(\ref{osc.theo.topo}) are trivially
implied by the single condition
$$
\int_{B_{\varrho_k}(x_k)} |A|^2 \d\mu < \ve^2 \quad \mbox{ for all }k=1,\ldots,K.
$$
In fact, the estimate of the conformal factor can then be shown in any codimension.

\setcounter{equation}{0}

\section{Estimation modulo the M\"obius group}
\label{moeb}

It will be essential in Theorem \ref{moeb.theo} 
to pass to a good representative under the action 
of the M\"obius group. The following lemma 
yields the desired M\"obius transformation.

\begin{lemma} \label{moeb.arrange}
Let $f: \sur \rightarrow \rel^n$ be an immersion of a closed surface $\sur$, 
with conformally invariant energy $\int_\sur |A^\circ|^2 \d \mu =: E$. Then 
there exists a M\"obius transformation $\phi$ such that $\tilde{f} = \phi \circ f$ 
satisfies $\tilde{f}(\sur) \subseteq B_1(0)$ and
\begin{equation}
\label{eqmoeb.arrange}
	\int_{B_{\varrho_0}(x)} |\tilde{A}^\circ|^2 \d \tilde{\mu} \leq E/2
	\quad \mbox{for all } x \in \rel^n, 
	\mbox{ where } \varrho_0 = \varrho_0(n,E) > 0.
\end{equation}
\end{lemma}
{\pr Proof:} 
By a dilation we may assume that for all $x \in \rel^n$ and some $x_0 \in \rel^n$ we have 
\begin{equation} \label{moeb.arrange.aux}
\int_{B_1(x)} |A^\circ|^2\,\d \mu \leq E/2 \leq
\int_{\overline{B_1(x_0)}} |A^\circ|^2\,\d \mu. 
\end{equation}
From (\ref{intro.gauss}) we see that the total Willmore energy of $f$ is bounded by 
\ee \label{moeb.arrange.curv-gauss}
{\cal W}(f) \leq E/2 + 4\pi.
\ef
We now prove by area comparison arguments that there is a point $x \in \rel^n$ satisfying
\ee \label{moeb.arrange.seek}
	\overline{B_1(x)} \cap f(\sur) = \emptyset
	\quad \mbox{and} \quad
	|x-x_0| \leq C(n,E).
\ef
The Li-Yau-inequality as in (\ref{pre.li-yau.dense}) yields the upper bound
\ee \label{moeb.arrange.dense-ratio}
	r^{-2}\mu(B_r(x_0)) \leq C(E)
	\quad \mbox{for any } r > 0,
\ef
while $\varrho = 1$ and $\sigma \searrow 0$ in (\ref{pre.li-yau.mono}) yields
\ee \label{moeb.arrange.lower}
\mu(B_1(x)) + {\cal W}(f,B_1(x)) \geq c > 0 \quad
\mbox{for any } x \in f(\sur).
\ef
For $R > 0$ to be chosen, let $B_2(x_j)$, $j=1,\ldots, N$, 
be a maximal disjoint collection of $2$-balls with centers
$x_j \in B_R(x_0)$. As the balls $B_4(x_j)$ cover $B_R(x_0)$ we have 
$N \geq R^n/4^n$. If $f(\sur) \cap \overline{B_1(x_j)} \neq \emptyset$ for all $j$, then
(\ref{moeb.arrange.lower}), (\ref{moeb.arrange.dense-ratio}) and (\ref{moeb.arrange.curv-gauss})
imply 
$$
c\, N \leq \sum_{j=1}^N \Big(\mu(B_2(x_j)) + {\cal W}(f,B_2(x_j))\Big)
\leq C(E)\,(R^2+1),
$$
thus $R \leq C(n,E)$. Taking $R = C(n,E) +1$ yields 
(\ref{moeb.arrange.seek}) for appropriate $x=x_j$.\\
\\
Translating by $-x$, we can assume that $x=0$ in (\ref{moeb.arrange.seek}), that is 
$f(\sur) \subseteq \rel^n - \overline{B_1(0)}$. For $R: = C(n,E)+1$ with $C(n,E)$ 
as in (\ref{moeb.arrange.seek}), we obtain from (\ref{moeb.arrange.aux}) for all 
$x \in \rel^n$
\ee \label{moeb.arrange.aux2}
\ad\D
\int_{B_1(x)} |A^\circ|^2 \d \mu \leq E / 2, \quad \mbox{ and } \quad
\int_{\rel^n - B_R(0)} |A^\circ|^2 \d \mu \leq E / 2.
\af
\ef
Now for $\tilde f = \phi \circ f$ where $\phi(x) = x/|x|^2$ we clearly have
$\tilde f(\sur) \subseteq B_1(0)$. Moreover if $|x| \geq 1/(2R)$, then 
a ball $B_\varrho(x)$ of radius $\varrho = \frac{1}{2}(\sqrt{1+ R^{-2}}-1)$ is mapped 
by $\phi^{-1} = \phi$ to a ball $B_{\varrho^\ast}(x^\ast)$ with 
$\varrho^\ast \leq 1$, and claim (\ref{eqmoeb.arrange}) follows
from (\ref{moeb.arrange.aux2}) using that the integral is locally conformally 
invariant. In the remaining case $|x| \leq 1/(2R)$,
we use $B_\varrho(x) \subseteq B_{1/R}(0)$ for $\varrho \leq 1/(2R)$, 
and obtain (\ref{eqmoeb.arrange}) from the second inequality 
in (\ref{moeb.arrange.aux2}). 
\proof
We can now prove our main theorem, recalling from (\ref{eqomegaconstants}) 
the definition of the $\omega^n_p$. 

\begin{theorem} \label{moeb.theo}
For $n =3,4$ and $p \geq 1$, let $\C(n,p,\delta)$ be the class of closed, 
oriented, genus $p$ surfaces $f: \sur \rightarrow \rel^n$ satisfying 
$\W(f) \leq \omega^n_p - \delta$ for some $\delta > 0$. 
Then for any  $f \in \C(n,p,\delta)$ there is a
M\"obius transformation $\phi$ and a constant curvature metric $g_0$,
such that the metric $g$ induced by $\phi \circ f$ satisfies
$$
g = e^{2u}g_0 \quad \mbox{ where }\quad
\max_{\sur} |u| \leq C(p,\delta) < \infty.
$$
\end{theorem}
{\pr Proof:}
It is obviously sufficient to obtain the result for small $\delta > 0$.
We consider an arbitrary sequence of surfaces $f_j \in \C(n,p,\delta)$, 
and put $g_j = f_j^{\ast}g_{euc}$, $\mu_j = f_j(\mu_{g_j})$. 
From (\ref{intro.gauss}) we have
\ee \label{moeb.theo.enrg-aux}
\int_\sur |A^\circ_j|^2 \d \mu_{g_j}
\leq 2(\omega^n_p - \delta) + 8 \pi (p - 1) \leq 8 \pi (p + 1) -2\delta.
\ef
Using Lemma \ref{moeb.arrange} we may assume after applying
suitable M\"obius transformations that
\ee \label{moeb.theo.normal}
\ad\D
f_j(\sur) \subseteq B_1(0) \quad \mbox{ and } \quad 
\int_{B_{\varrho_0}(x)} |A^\circ_j|^2 \d \mu_j
\leq 4 \pi (p + 1) - \delta \quad \mbox{ for all }x \in \rel^n,
\af
\ef
where $\varrho_0 > 0$ depends only on the genus $p$. The uniformization 
theorem yields unique conformal metrics $e^{-2u_j} g_j$ having the same 
area and constant curvature. The theorem will be proved by showing that 
\ee \label{moeb.theo.stat}
	\liminf_{j \rightarrow \infty} \|u_j\|_{L^\infty(\sur)} < \infty.
\ef
We start by recalling from  (\ref{pre.li-yau.dense}) the Li-Yau-inequality
\ee \label{moeb.theo.dense-ratio}
	\varrho^{-2} \mu_j(B_\varrho(x)) \leq C
	\quad \mbox{for all } x \in \rel^n,\,\varrho > 0.
\ef
For $\alpha_j = f_j(\mu_{g_j} \llcorner |A_j|^2)$ we have $\alpha_j(\rel^n) \leq C(p)$,
hence for a subsequence
\ee \label{moeb.theo.conv}
	\mu_j, \alpha_j \rightarrow \mu, \alpha
	\quad \mbox{weakly}^* \mbox{ in } C_c^0(\rel^n)^*.
\ef
Moreover we see as in \bcite{sim.will} p. 310 that
\ee \label{moeb.theo.haus}
    \spt \mu_j \rightarrow \spt \mu
    \quad \mbox{in Hausdorff distance},
\ef
which yields further 
$\spt \alpha \subseteq \spt \mu \subseteq \overline{B_1(0)}$.
Now by Allard's integral compactness theorem for varifolds, see
\bcite{sim} Remark 42.8, the measure $\mu$ is an integral $2$-varifold
with weak mean curvature $H_\mu \in L^2(\mu)$, more precisely we have
$$
\W(\mu):= \frac{1}{4} \int |H_\mu|^2 \d\mu \leq \liminf_{j \to \infty} \W(f_j) \leq 8\pi-\delta.
$$
As discussed in the appendix of \cite{kuw.schae.will3}, the monotonicity 
formula from \cite{sim.will} applies to varifolds with weak mean curvature in $L^2$,
in particular the Li-Yau inequality (\ref{pre.li-yau}) yields 
\ee \label{moeb.theo.dense}
\theta^2(\mu,x) \leq \frac{8 \pi - \delta}{4 \pi} = 2 - \frac{\delta}{4 \pi}
\quad \mbox{for all }x \in \rel^n.
\ef
We further obtain, writing $\perp$ for the projection onto $(T_x\mu)^\perp$,
\ee \label{pre.li-yau.pos}
\int_{B_\sigma(x_0)} \frac{|(x-x_0)^\perp|^2}{|x-x_0|^4}\d \mu(x) < \infty
\quad \mbox{ for all } x_0 \in \rel^n.
\ef
Let $\varepsilon_0 = \varepsilon_0(n,\beta)$ be the constant from Lemma 2.1 of 
\bcite{sim.will}; we take $\beta = C$ for $C > 0$ as in (\ref{moeb.theo.dense-ratio}) 
whence $\varepsilon_0 > 0$ is universal. For $\varepsilon_1 \in  (0,\varepsilon_0]$ there
are only finitely many points $x_1, \ldots, x_K$ with
\dd
\alpha(\{x_k\}) \geq \varepsilon_1^2
\quad \mbox{for } k = 1, \ldots, K,
\df
in fact $K \leq C(p) \varepsilon_1^{-2}$. For given
$\varepsilon \in (0,\varepsilon_1)$ we may use (\ref{pre.li-yau.pos})
and (\ref{moeb.theo.dense}) to choose 
$\varrho \in (0,\frac{1}{2}\min_{k \neq l} |x_k - x_l|)$ with $\varrho < \varrho_0$,
such that for all $k$ we have the inequalities
\dcd
\alpha(\overline{B_\varrho(x_k)} - \{x_k\}) < \varepsilon^2, \\
\mu(\overline{B_{7 \varrho / 8}(x_k)})
< (2 - \frac{\delta}{20}) \pi (7 \varrho / 8)^2, \\
\displaystyle{\int_{\overline{B_\varrho(x_k)}}
\frac{|(x-x_k)^\perp|^2}{|x-x_k|^4}  \d \mu(x) < \varepsilon^2}.
\dcf
For any $y \notin \{x_1,\ldots,x_K\}$ 
there exists a radius $\varrho_y \in (0,\varrho_0)$ 
such that $\alpha(\overline{B_{\varrho_y}(y)}) < \varepsilon_1^2$.
Now we select finitely many points $y_1, \ldots, y_L \in \spt \mu
- \bigcup_{k=1}^K B_{\varrho/2}(x_k)$ such that
\dd
\spt \mu \subseteq
\bigcup_{k=1}^K B_{\varrho/2}(x_k)
\cup \bigcup_{l=1}^L B_{\varrho_{y_l}/2}(y_l).
\df
By (\ref{moeb.theo.conv}) and (\ref{moeb.theo.haus})
we get for any $r \in (0,\varrho/2]$ and
$j$ sufficiently large (depending on $r$)
\ee \label{moeb.theo.aux}
\left.
\ad
	\D f_j(\sur) = \spt \mu_j \subseteq 	\bigcup_{k=1}^K B_{\varrho/2}(x_k) \cup \bigcup_{l=1}^L B_{\varrho_{y_l}/2}(y_l), \\
	\D \int_{B_\varrho(x_k) - B_r(x_k)} |A_j|^2 \d \mu_j < \varepsilon^2, \\
	\D \int_{B_{\varrho_{y_l}}(y_l)} 	|A_j|^2 \d \mu_j < \varepsilon_1^2, \\
	\D \mu_j(B_{7 \varrho / 8}(x_k)) < (2 - \delta / 20) \pi (7 \varrho / 8)^2, \\
	\D \int_{B_\varrho(x_k) - B_r(x_k)} {\frac{|(x-x_k)^\perp|^2}{|x-x_k|^4}} \d \mu_j(x) < \varepsilon^2, 
\af
\right\}
\ef
for $k = 1, \ldots, K$ and $l = 1, \ldots, L$. For the covering 
in (\ref{moeb.theo.aux}) we shall now verify the assumptions
of Theorem \ref{osc.theo}, provided
that $\varepsilon_1 = \varepsilon_1(n,\delta) \in (0,\varepsilon_0)$ and 
$\varepsilon = \varepsilon(n,p,\delta) \in (0,\varepsilon_1)$ 
are sufficiently small.\\ 
\\
The condition (\ref{osc.theo.topo}) is clearly 
satisfied on the $B_{\varrho_{y_l}}(y_l)$, $l = 1, \ldots, L$,
for $\varepsilon_1 > 0$ sufficiently small. For $k \in \{1, \ldots, K\}$,
we have the assumptions of Lemma \ref{lemma.pre.graph} and also
(\ref{pre.graph.densityass}) for any $r \in (0,\varrho/2]$. Thus
for the multiplicity $M_k$ as in (\ref{pre.graph.mult}), we get
from (\ref{pre.graph.density})
\begin{eqnarray*}
(1 - C \varepsilon) M_k \pi \Big( (7 \varrho / 8)^2 - (5 r / 4)^2 \Big)
& \leq & \mu_j(B_{7 \varrho /8}(x_k) - B_{5 r / 4}(x_k))\\ 
& \leq & \mu_j(B_{7 \varrho / 8}(x_k))\\
& \leq & (2 - \delta / 20) \pi (7 \varrho / 8)^2,
\end{eqnarray*}
Assuming $\varepsilon \leq \varepsilon(\delta)$ and $r / \varrho \leq c(\delta)$
this implies 
\ee \label{moeb.theo.bound}
	M_k = 1.
\ef
For $\sigma \in [5\varrho_k/8,7\varrho_k/8]$ as in Lemma \ref{lemma.pre.graph},  
we conclude that $f_j^{-1}\big(B_\sigma(x_k)\big)$ is bounded by just one 
circle, and can be compactified to a closed surface $\sur_k = \sur_{j,k}$ 
of genus $p_k = p_{j,k}$ by adding one disc. This means we have
\ee \label{moeb.theo.genus-k}
\chi(f_j^{-1}(B_\sigma(x_k))) = 2(1-p_k) -1.
\ef
As (\ref{eqgeodesiccurvaturebound}) holds with multiplicity one, the 
Gau{\ss}-Bonnet theorem yields
\begin{equation}
\label{moeb.theo.gaussintegral}
\Big| \int_{B_\sigma(x_k)} K_j \d\mu_j + 4\pi p_k\Big| \leq C \ve^\alpha.
\end{equation}
Now $K_j \geq \frac{1}{2} |A^\circ_j|^2$ by (\ref{intro.gauss1}), and using 
$\varrho \leq \varrho_0$ we see from (\ref{moeb.theo.enrg-aux}) and 
(\ref{moeb.theo.normal}) that
$$
\int_{B_\sigma(x_k)} K_j \d\mu_j \geq 
-\frac{1}{2} \int_{B_{\varrho_k}(x_k)} |A^\circ_j|^2 \d\mu_j \geq 
-\frac{1}{4}\int_{\sur} |A^\circ_j|^2 \d\mu_j \geq
-2\pi(p+1) + \frac{\delta}{2}.
$$
In the case $p=1$ this implies the condition (\ref{osc.theo.topo}) with $\delta/2$ 
instead of $\delta$ as well as $p_k = 0$, for $\varepsilon \leq \varepsilon(n,\delta)$.
For $p \geq 2$ we get 
\ee \label{moeb.theo.genus-top}
p_k < (p+1)/2 < p
\quad \mbox{for } k = 1, \ldots, K.
\ef
For appropriate $\sigma_k = \sigma_{j,k} \in ]5 \varrho/8,7 \varrho/8[$, 
we now use Lemma \ref{lemma.pre.graph}(b) to attach an end to the 
restriction of $f_j$ to $f_j^{-1}(B_{\sigma_k}(x_k))$, obtaining
an immersion $\tilde{f}_{j,k}: \sur_k - \{q_k\} \rightarrow \rel^n$
such that under $\tilde{f}_{j,k}$ a neighborhood of the puncture $q_k$ 
corresponds to a neighborhood of infinity in some affine plane, and such that 
\ee \label{moeb.theo.ext}
\int_{\rel^n - B_{\sigma_k}(x_k)} |\tilde{A}_k|^2 \d \tilde{\mu}_k
\leq C \varepsilon^2.
\ef
By (\ref{moeb.theo.ext}) and the conformal invariance of the Willmore energy, we get
\ee \label{moeb.theo.ext-ener}
{\cal W}(f_j , B_{\sigma_k}(x_k))
\geq {\cal W}(\tilde{f}_{j,k}) - C \varepsilon^2
\geq \beta_{p_k}^n - 4 \pi - C \varepsilon^2.
\ef
Adding $k$ discs to $\sur - \bigcup_{k=1}^K B_{\sigma_k}(x_k)$ yields a surface
of some genus $p_0$, where 
\begin{eqnarray*}
2 (1 - p) = \chi(\sur) & = &
\chi \Big(\sur - \bigcup_{k=1}^K f_j^{-1}(B_{\sigma_k}(x_k)) \Big)
+ \sum_{k=1}^K \chi \left(f_j^{-1}(B_{\sigma_k}(x_k)) \right)\\
& = &  2 (1 - p_0) - K + \sum_{k=1}^K \Big(2 (1 - p_k) - 1\Big)\\
& = & 2 \Big(1-\sum_{k=0}^K p_k\Big),
\end{eqnarray*}
which means 
\begin{equation}
\label{eqpartition}
p = \sum_{k=0}^K p_k.
\end{equation}
In fact, adding the discs with bounds as in (\ref{moeb.theo.ext}), we see that 
\ee \label{moeb.theo.ext-ener0}
{\cal W}\Big(f_j ,\rel^n - \bigcup_{k=1}^K B_{\sigma_k}(x_k) \Big)
\geq \beta_{p_0}^n - C(K) \varepsilon^2.
\ef
Combining (\ref{moeb.theo.ext-ener}) and (\ref{moeb.theo.ext-ener0}) implies
$$
\sum_{k=0}^K (\beta_{p_k}^n - 4 \pi) \leq 
{\cal W}(f_j) - 4\pi + C(K) \varepsilon^2
\leq \omega^n_p - \delta - 4\pi + C (K) \varepsilon^2
< \omega_p^n - 4\pi,
$$
if $\varepsilon \leq \varepsilon(K,\delta)$. From (\ref{eqpartition}), (\ref{moeb.theo.genus-top})
and the definition of the $\omega^n_p$, see (\ref{eqomegaconstants}), we now see
\begin{equation}
\label{eqtopology} 
p_0 = p, \quad \mbox{ and } \quad  p_k = 0 \quad \mbox{for }k=1,\ldots,K.
\end{equation}
Together with (\ref{moeb.theo.gaussintegral}) this establishes (\ref{osc.theo.topo}) for
any $p \geq 1$.\\
\\
Next, claim (\ref{osc.theo.curv-ann}) is immediate by taking 
$\varepsilon,\varepsilon_1 \leq \varepsilon(\Lambda,\delta)$ and 
$r \leq \varrho/2$ in (\ref{moeb.theo.aux}). Moreover,
for $l = 1, \ldots, L$ we get from (\ref{moeb.theo.aux})
\dd
2 \int_{B_{\varrho_{y_l}}(y_l)} |K_j| \d \mu_j,\,
\int_{B_{\varrho_ {y_l}}(y_l)} |A_j^\circ|^2 \d \mu_j
\leq \int_{B_{\varrho_{y_l}}(y_l)} |A_j|^2 \d \mu_j
\leq \varepsilon_1^2,
\df
hence (\ref{osc.theo.curv3}) and (\ref{osc.theo.curv4}) hold
for $\varepsilon_1 > 0$ small enough.
For $k = 1, \ldots, K$ we get from (\ref{eqgeodesiccurvaturebound})
combined with (\ref{moeb.theo.bound}), (\ref{moeb.theo.genus-k})
and (\ref{eqtopology}), for appropriate
$\sigma \in ]5 \varrho / 8 , 7 \varrho / 8[$,
\begin{equation}
\label{eqgaussestimate}
\Big| \int_{B_\sigma(x_k)} K_j \d \mu_j\Big|
\leq 2\pi \left| \chi\Big(f_j^{-1}(B_\sigma(x_k))\Big)-1\right| + C \varepsilon^\alpha
	= C \varepsilon^\alpha.
\end{equation}
From $|K| \leq \frac{1}{2}|A|^2 = |A^\circ|^2 + K$ we have the inequality
$$
\int_{B_\varrho(x_k)} |K_j| \d \mu_j \leq 
\frac{1}{2} \int_{B_\varrho(x_k) - B_{\varrho/2}(x_k)} |A_j|^2 \d \mu_j
+ \int_{B_\sigma(x_k)} |A_j^\circ|^2 \d \mu_j
+ \int_{B_\sigma(x_k)} K_j \d \mu_j,\\
$$
hence we obtain from (\ref{moeb.theo.aux}) and (\ref{eqgaussestimate})
\ee 
\label{moeb.theo.curv.aux}
\int_{B_\varrho(x_k)} |K_j| \d \mu_j \leq
\int_{B_\varrho(x_k)} |A_j^\circ|^2 \d \mu_j  + C \varepsilon^\alpha.
\ef
We proceed similarly using (\ref{moeb.theo.aux}), (\ref{eqgaussestimate}),
(\ref{moeb.theo.ext-ener0}), (\ref{eqtopology}) and 
$|A^\circ|^2 = |\mean|^2/2 - 2K$  
\begin{eqnarray*}
\frac{1}{2} \int_{B_\varrho(x_k)} |A_j^\circ|^2 \d \mu_j & \leq &
\frac{1}{4} \int_{B_\sigma(x_k)} |\mean_j|^2 \d \mu_j + C \varepsilon^\alpha\\
& \leq & {\cal W}(f_j) - {\cal W}(f_j,\rel^n - \bigcup_{k=1}^K B_{\sigma_k}(x_k)) + C \varepsilon^\alpha\\
& \leq & \omega^n_p - \delta - \beta^n_p + C \varepsilon^\alpha.
\end{eqnarray*}
As $\omega^n_p \leq 8\pi$ and $\beta^n_p \geq 4\pi$, we conclude
from (\ref{moeb.theo.curv.aux})
$$
\int_{B_\varrho(x_k)} |K_j| \d \mu_j 
\leq 2( \omega^n_p - \beta^n_p) - 2\delta + C \varepsilon^\alpha 
\leq 8\pi - 2\delta + C\varepsilon^\alpha,
$$
which proves (\ref{osc.theo.curv3}) taking $\varepsilon \leq \varepsilon(\delta)$.
For $n=4$ we have 
\begin{eqnarray*}
\int_{B_\varrho(x_k)} |K_j| \d \mu_j
+ \frac{1}{2} \int_{B_\varrho(x_k)}  |A_j^\circ|^2 \d \mu_j 
& \leq &\frac{3}{2} \int_{B_\varrho(x_k)} |A_j^\circ|^2 \d \mu_j + C \varepsilon^\alpha\\
& \leq & 3(\omega^4_p - \delta - \beta_p^4) + C \varepsilon^\alpha. 
\end{eqnarray*}
Now (\ref{osc.theo.curv4}) follows by definition of $\omega^4_p$ for 
$\varepsilon \leq \varepsilon(\delta)$ small enough. Thus all conditions of 
Theorem \ref{osc.theo} are verified, and application of that theorem 
finishes the proof.
\proof

\setcounter{equation}{0}
\section{Compactness in moduli space} \label{modul}

The main result of this section is 

\begin{theorem} For $n \in \{3,4\}$ and $p \geq 1$,
the conformal structures induced by immersions $f$ in $\C(n,p,\delta)$
are contained in a compact subset $K = K(p,\delta)$ of the moduli space.
\end{theorem}
The theorem follows directly from Theorem \ref{moeb.theo} and the following 

\begin{lemma} \label{lemmamoduli} Let $f:\sur \to \rel^n$ be an immmersion
of a closed oriented surface of genus $p \geq 1$, with induced metric 
$g = f^\ast g_{euc}$. Assume that
$$
\W(f),\, \max_{\sur} |u| \leq \Lambda,
$$
where $\W(f)$ is the Willmore energy and $g_0: = e^{-2u}g$ is a conformal
metric of constant curvature. Then the conformal structure induced by $g$ 
lies in a compact subset $K = K(n,p,\Lambda)$ of the moduli space.
\end{lemma}
{\pr Proof: } We first give the proof for $p \geq 2$, where we normalize to 
$K_{g_0} \equiv -1$ by a dilation. 
Let $\ell > 0$ be the length of a shortest closed geodesic in $(\sur,g_0)$.
By the Mumford compactness theorem, see e.g. \cite{tromba} Theorem C.1, the
lemma follows from a lower bound for $\ell$ depending only on $n,p$ 
and $\Lambda$. As the hyperbolic plane has no conjugate points, 
we have ${\rm inj}(M,g_0) = \ell/2$ by an argument of Klingenberg, 
see Lemma 4 in \cite{klinge}, and hyperbolic geometry implies
\begin{equation}
\label{eqareaballs}
\mu_{g_0}(B^{g_0}_r(p)) \geq \pi r^2 \quad \mbox{ for all }0 < r \leq \ell/2.
\end{equation}
Select a closed geodesic $\gamma$ for $g_0$ of length $\ell$. With respect
to geodesic distance, there is a parallel neighborhood of $\gamma$
which is isometric to the quotient
of $\{re^{i\theta}: r > 0,\,|\theta-\pi/2| < \theta_0\}$ by the action 
of $e^{k\ell}$, $k \in {\mathbb Z}$, where $\gamma$ corresponds to $\theta = \pi/2$. 
Clearly $\gamma$ is not contractible since otherwise it would lift 
to a closed geodesic in the hyperbolic plane.
By the collar lemma, see \cite{tromba} Lemma D.1, we may take
$\theta_0 \in (0,\pi/4]$ as a universal constant, as we can assume 
without loss of generality that $\ell \leq 1$. Now let 
$p_1 \simeq e^{i\theta_1},\ldots,p_K \simeq e^{i\theta_K}$ be a maximal collection 
of points with $|\theta_j - \pi/2| < \theta_0$, such that 
the balls $B^{g_0}_{\ell}(p_j)$ are pairwise disjoint. By maximality the 
$B^{g_0}_{2\ell}(p_j)$ cover the geodesic $\{e^{i\theta}: |\theta-\pi/2| < \theta_0\}$, 
which implies that $K \geq c_0/\ell$ for a universal constant $c_0 > 0$. 
The closed curves $\gamma_k$ corresponding to $e^{t+i\theta_k}$, $0 \leq t \leq \ell$, 
have length $L_{g_0}(\gamma_k) \leq C\ell$. We conclude
\begin{equation}
\label{eqlengthestimate}
L_g(\gamma_k) \leq C(\Lambda) L_{g_0}(\gamma_k) \leq C(\Lambda) \ell =:\varrho/4.
\end{equation}
Given $k \in \{1,\ldots,K\}$, we denote by $I_k$ the set of those $i \in \{1,\ldots,K\}$ 
for which $f(p_i) \in B_{2\varrho}(f(p_k))\}$.  For $i \in I_k$ and 
$p \in B^{g_0}_\ell(p_i)$ we estimate
$$
|f(p)-f(p_k)| \leq {\rm dist}_g(p,p_i) + |f(p_i)-f(p_k)|
\leq C(\Lambda) {\rm dist}_{g_0}(p,p_i) + 2\varrho
\leq C(\Lambda) \varrho.
$$
As the balls $B^{g_0}_\ell(p_i)$ are pairwise disjoint, we get 
putting $r = \ell/2$ in (\ref{eqareaballs})
$$
(\# I_k)\frac{\pi \ell^2}{4} \leq \sum_{i \in I_k} \mu_{g_0}(B^{g_0}_\ell(p_i))
\leq C(\Lambda)\, \mu_g\Big(f^{-1}B_{C(\Lambda)\varrho}\big(f(p_k)\big)\Big)
\leq C(\Lambda)\varrho^2,
$$
where the last step uses the Li-Yau inequality (\ref{pre.li-yau.dense}). We thus have 
\begin{equation}
\label{eqcoverestimate}
\# I_k \leq C(\Lambda) \quad \mbox{ for }k=1,\ldots,K.
\end{equation}
Now choose a maximal set $J \subseteq \{1,\ldots,K\}$ with 
$B_\varrho(f(p_k)) \cap B_\varrho(f(p_l)) = \emptyset$ for $k \neq l$.
For any $m \in \{1,\ldots,K\}$ we have $f(p_m) \in B_{2\varrho}(f(p_k))$ 
for some $k \in J$, which means $\{1,\ldots,K\} = \bigcup_{k \in J} I_k$.
By (\ref{eqcoverestimate}) this yields
$K \leq \sum_{k \in J} \# I_k \leq C(\Lambda) \# J$ and hence
\begin{equation}
\label{eqinjectivityestimate}
\#J \geq c_0/\ell \quad \mbox{ for }c_0 = c_0(\Lambda) > 0.
\end{equation}
As the $B_\varrho(f(p_k))$ are disjoint for $k \in J$, we get for some $k \in J$
using Gau{\ss}-Bonnet
$$
\int_{B_\varrho(f(p_k))} |A|^2 \d\mu \leq \frac{1}{\# J} \int_{\sur} |A|^2\,\d\mu
\leq C(\Lambda,p)\ell.
$$
Thus for $C(\Lambda,p)\ell < \ve_0(n,\Lambda)$ the assumptions of \cite{sim.will} 
Lemma 2.1 are satisfied, recalling also the density ratio estimate 
(\ref{pre.li-yau.dense}), hence there exists a $\sigma \in ]\varrho/4,\varrho/2[$ 
such that $f^{-1}\big(B_\sigma(f(p_k))\big)$ is a disjoint union of discs 
$D^i_\sigma$, $i=1,\ldots,M$. Now by (\ref{eqlengthestimate}) we have 
$L_g(\gamma_k) \leq \varrho/4 < \sigma$ which implies that $f \circ \gamma_k$ 
lies in $B_\sigma((f(p_k))$, or equivalently $\gamma_k$ is contained in 
$f^{-1}\big(B_\sigma(f(p_k))\big)$. But then $\gamma_k$ is actually contained 
in one of the discs $D^i_\sigma$, in particular $\gamma_k$ is contractible 
in $\sur$. But then $\gamma$ is also contractible which contradicts our 
previous observation.\\
\\
For $p =1$ we normalize such that $\mu_{g_0}(\sur) = 1$. It is well-known that
$(\sur,g_0)$ is isometric to the quotient of $\rel^2$ by a lattice of the form
$\Gamma/\sqrt{b}$, where $\Gamma = \ganz + \ganz (a,b)$ with
$0 \leq a \leq 1/2$, $a^2 + b^2 \geq 1$ and $b > 0$; here dilating the lattice 
by $1/\sqrt{b}$ adjusts the volume to one. The length of a shortest closed
geodesic is then $\ell = 1/\sqrt{b}$, in fact any horizontal line segment of
that length corresponds to a shortest closed geodesic. 
We now consider points $p_k$
corresponding to $(0,2k\ell)$ for $k=1,\ldots,K$. It is elementary that 
we can achieve $B^{g_0}_\ell(p_k) \cap B^{g_0}_\ell(p_l) = \emptyset$ 
where $K \geq c_0 \ell^{-2}$. The horizontal segments yield closed 
geodesics $\gamma_k$ through $p_k$ of length $L_{g_0}(\gamma_k) = \ell$. From 
here the proof proceeds as in the case $p \geq 2$.
\proof
We finally discuss the optimality of the constants $\omega^n_p$ in 
Theorem \ref{moeb.theo}. A standard example, see \cite{sim.will}, 
is obtained by connecting two concentric round spheres at small distance 
by $p+1$ suitably scaled catenoids. This yields a sequence of 
embeddings $f_j:\sur \to \rel^3$ of genus $p \geq 1$ with $\W(f_j) \to 8\pi$.
By a dilation we have in addition that $\mu_{g_j}(\sur) = 1$ for all $j$.
Assume by contradiction that there exist M\"obius transformations 
$\phi_j$ and constant curvature metrics $g_{0,j}$, such that 
$\tilde{g}_j = (\phi_j \circ f_j)^\ast g_{euc} = e^{2u_j}g_{0,j}$ where 
$\max_{\sur}|u_j|$ remain in a compact set as $j \to \infty$. Composing $\phi_j$ 
with a suitable dilation we may assume that $\mu_{\tilde{g}_j}(\sur) = 1$. 
By Lemma \ref{lemmamoduli}, the conformal structures induced by 
the $g_j$ remain bounded, which implies that the minimal length 
of a noncontractible loop with respect to $g_{0,j}$, and hence 
with respect to $\tilde{g}_j$, is bounded below independent of $j$.
In particular, the metric $\tilde{g}_j$ is not uniformly bounded 
by $g_j$ near the concentrating catenoids.
Now $\phi_j$ is a composition of a Euclidean motion,
a dilation and an inversion, hence we have 
$$
\tilde{g}_j = c_j^2\, g_j \quad \mbox{ or }\quad
\tilde{g}_j = \frac{c_j^2}{|f_j-a_j|^4}\, g_j, \quad
\mbox{ where }c_j > 0,\,a_j \in \rel^3.
$$
In the first case, the area normalization implies $c_j = 1$ which 
is a contradiction. In the second case, we note that the 
$a_j$ cannot diverge since otherwise we get for large $j$
$$
\frac{c_j^2}{16 |a_j|^4}\, g_j \leq \tilde{g}_j \leq \frac{16 c_j^2}{|a_j|^4}\,g_j.
$$
The area normalization yields $1/16 \leq c_j^2/|a_j|^4 \leq 16$, and we 
have a contradiction as before. Thus we can assume that the $a_j$ converge
to some $a \in \rel^3$, and also that the $c_j$ remain bounded. But since
$p+1 \geq 2$ there is a catenoid concentrating at a point different to $a$, 
and at that point $\tilde{g}_j$ remains bounded by $g_j$. This contradiction 
shows that the constant $\omega^n_p$ in Theorem \ref{moeb.theo} cannot be replaced 
by a constant strictly bigger than $8\pi$. Inverting surfaces of genus $p_i$
where $p_1+\ldots+p_k = p$ at points on the surface and then glueing them 
into a round sphere, we see similarly that $\omega^n_p$ cannot be replaced by 
a number bigger than $\tilde{\beta}^n_p$, 
and in particular that $\omega^3_p = \min\{8\pi,\beta^3_p\}$ is optimal 
for the statement of Theorem \ref{moeb.theo}.

%
\setcounter{equation}{0}
\section{Conformal parametrization} \label{conf}

In this section, we prove the estimate for the conformal factor needed in 
the proof of Theorem \ref{osc.theo}, thereby extending results 
of \bcite{muell.sver}.

\begin{theorem} \label{conf.flat}

Let $f: \rel^2 \rightarrow \rel^n$, $n = 3,4$, be a complete 
conformal immersion with induced metric $g = e^{2u} g_{euc}$
and square integrable second fundamental form satisfying
\ee \label{conf.flat.curv}
\int_{\rel^2}  K \d \mu_g = 0 \quad \mbox{ for }K = K_g,
\ef
\ee \label{conf.flat.curv3}
\int_{\rel^2}  |K| \d \mu_g \leq 8 \pi - \delta \quad \mbox{for }n=3,
\ef
\ee \label{conf.flat.curv4}
\left.
\ad\D
\int_{\rel^2}  |K| \d \mu_g + 
\frac{1}{2} \int_{\rel^2} |A^\circ|^2 \d \mu_g \leq 8 \pi - \delta, \\
\int_{\rel^2} |A^\circ|^2 \d \mu_g < 8 \pi,
\af
\right\}
	\quad \mbox{for } n = 4,
\ef
for some $\delta > 0$. Then the limit $\lambda = \lim_{z \rightarrow \infty} u(z) \in \rel$ 
exists, and 
\ee \label{conf.flat.esti}
\|u - \lambda\|_{L^\infty(\rel^2)},
\|D u\|_{L^2(\rel^2)},
\|D^2 u\|_{L^1(\rel^2)}
\leq C(\delta) \int_{\rel^2} |A|^2 \d \mu_g.
\ef 
\defin
\end{theorem}
We shall prove this theorem by constructing a solution 
$v: \rel^2 \rightarrow \rel$ of the problem
\ee \label{conf.flat-aux.equ}
- \Delta_g v = K  \quad \mbox{in } \rel^2, \quad \mbox{ and } \quad
\lim_{z \rightarrow \infty} v(z) = 0,
\ef
which satisfies the estimates 
\ee \label{conf.flat-aux.esti}
\|v \|_{L^\infty(\rel^2)},
\|D v\|_{L^2(\rel^2)},
\|D^2 v\|_{L^1(\rel^2)}
\leq C(\delta) \int_{\rel^2} |A|^2 \d \mu_g.
\ef 
The claim then follows easily. In fact, the function $u$ solves 
$-\Delta_g u = K_g$, see (\ref{eqconformalfactor}), 
hence the difference $u-v$ is an entire harmonic function. 
But \bcite{muell.sver} Theorem 4.2.1, Corollary 4.2.5,
combined with (\ref{conf.flat.curv}), imply that $u$ is bounded.
Therefore $u - v$ is also bounded and reduces to a constant $\lambda$,
which proves the lemma.\\
\\
{\pr Proof of Theorem \ref{conf.flat} for n = 3:}
The projection $\pi:S^3 \to \com P^1$, $(z_1,z_2) \mapsto [z_1:z_2]$, 
is a Riemannian submersion for the Fubini-Study metric $g_{FS}$ on $\com P^1$. 
Introduce the diffeomorphism $\cor: S^2 \to \com P^1$ induced by 
composing the standard chart 
$$
\psi:\com \to \com P^1,\, \psi(z) = \pi \Big(\frac{(z,1)}{\sqrt{|z|^2 +1}}\Big)
$$
with the stereographic projection
$$
T:S^2 \backslash \{-e_3\} \to \com,\,T(\zeta,s) = \frac{\zeta}{1+s}.
$$
One computes $\psi^\ast g_{FS} = (1+|z|^2)^{-2}g_{euc} = \frac{1}{4} (T^{-1})^\ast g_{S^2}$,
which implies 
\ee \label{conf.flat3.cor}
\cor^{\ast} g_{FS} = \frac{1}{4}\, g_{S^2}.
\ef
As the Jacobian of the normal $\nu: (\rel^2,g) \rightarrow S^2$ along $f$
is $J \nu = |K|$, we get by (\ref{conf.flat.curv3})
\dd
\int_{\rel^2} J(\cor \circ \nu) \d \mu_g
= \frac{1}{4} \int_{\rel^2} |K| \d \mu_g
\leq 2 \pi - \delta / 4.
\df
Recalling that the K\"ahler form $\omega$ on $\com P^1$, as defined in 
\bcite{muell.sver} 2.2, equals twice the volume form $vol_{FS}$ 
with respect to the Fubini-Study metric, 
we get $\cor^* \omega = 2 \cor^* vol_{FS} = \frac{1}{2} vol_{S^2}$ by (\ref{conf.flat3.cor}).
Hence using $\nu^* vol_{S^2} = K vol_g = K e^{2u} vol_{euc}$ we obtain
from (\ref{conf.flat.curv}) that
\dd
\int_{\rel^2} (\cor \circ \nu)^* \omega
= \int_{\rel^2} (K/2) vol_g = 0.
\df
We may therefore apply \bcite{muell.sver} Corollary 3.5.7 to get a solution
$v: \rel^2 \rightarrow \rel$ of 
$$
- \Delta v = *2 (\cor \circ \nu)^* \omega = K e^{2u} \quad  \mbox{on } \rel^2, \quad
\mbox{ with }\quad \lim_{z \to \infty} v(z) = 0,
$$
where $\Delta$, $\ast$ are taken with respect to the standard metric on $\rel^2$,
and such that 
\begin{eqnarray*}
\|v \|_{L^\infty(\rel^2)},\, \|D v\|_{L^2(\rel^2)},\, \|D^2 v \|_{L^1(\rel^2)}
& \leq &  C(\delta) \int_{\rel^2}  |D(\cor \circ \nu)|^2 \d \mu_g\\
& = & \frac{C(\delta)}{4} \int_{\rel^2}  |D \nu|^2 \d \mu_g\\
& = & \frac{C(\delta)}{4} \int_{\rel^2}  |A|^2 \d \mu_g.
\end{eqnarray*}
As $- \Delta_g v = K$ by construction, the lemma follows for $n=3$.
\proof

{\flushleft We }remark that if we use instead of $\cor$ the map $\tilde \cor$ identifying 
$S^2$ with the Gra{\ss}mannian $G_{3,2} \subseteq \com P^2$, then we have
$\tilde \cor^* g_{FS} = g_{S^2}/2$ instead of (\ref{conf.flat3.cor}), 
which implies only
\dd
\int_{\rel^2} J(\tilde \cor \circ \nu) \d \mu_g
= \frac{1}{2} \int_{\rel^2} |K| \d \mu_g,
\df
so that instead of (\ref{conf.flat.curv3}) we would need the stronger assumption
\dd
\int_{\rel^2} |K| \d \mu_g \leq 4 \pi - \delta  \quad \mbox{for } n = 3.
\df
For $n \geq 4$ the Jacobian $JG$ of the Gau{\ss} map 
$G: (\sur,g) \rightarrow G_{n,2} \subseteq \com P^{n-1}$
can in general not be expressed in terms of the Gau{\ss} curvature $K$ alone,
more precisely it was computed in \bcite{hoff.oss} that in points
where $\mean$ is nonzero one has
\dd
	JG = \frac{1}{2} \sqrt{|K|^2 + \frac{1}{2} |\mean|^2 |B|^2},
\df
where $B$ is the component of $A$ orthogonal to $\mean$. 
For the proof of Theorem \ref{conf.flat} for $n = 4$,
we will use a correspondance $G_{4,2} \leftrightarrow S^2 \times S^2$.
Recall that an oriented $2$-plane $P$ in $\rel^n$ with oriented orthonormal
basis $v,w$ is represented by 
\ee \label{conf.grass-def}
	[(v + i w)/\sqrt{2}]
	\in G_{4,2} = \Big\{[z_0: \ldots :z_{n-1}]\, \Big|
	\sum_{k=0}^{n-1} z_k^2 = 0\Big\}
	\subseteq \com P^{n-1}.
\ef
Alternatively, we can assign to $P$ the $2$-vector 
$v \wedge w \in \Lambda_2(\rel^n)$.
For $n = 4$ the Hodge operator $*: \Lambda_2(\rel^4) \to  \Lambda_2(\rel^4)$
is an involution, that is $*^2 = {\rm Id}$, and we have a direct sum decomposition
$\Lambda_2(\rel^4) = E_{+} \oplus E_{-}$ into the $\pm 1$ eigenspaces, with 
corresponding projections $\Pi_{\pm} \xi = (\xi \pm *\xi)/2$. As the Hodge star
is an isometry the decomposition is orthogonal, and both spaces $E_{\pm}$ 
are three-dimensional with orthonormal bases
\dcd
	e_{12}^+ := (e_1 \wedge e_2 + e_3 \wedge e_4) / \sqrt{2}, \quad
	e_{12}^- := (e_1 \wedge e_2 - e_3 \wedge e_4) / \sqrt{2}, \\
	e_{13}^+ := (e_1 \wedge e_3 + e_4 \wedge e_2) / \sqrt{2}, \quad
	e_{13}^- := (e_1 \wedge e_3 - e_4 \wedge e_2) / \sqrt{2}, \\
	e_{14}^+ := (e_1 \wedge e_4 + e_2 \wedge e_3) / \sqrt{2}, \quad
	e_{14}^- := (e_1 \wedge e_4 - e_2 \wedge e_3) / \sqrt{2}. \\
\dcf
We orient the $2$-spheres $S^2_{\pm} = S^5 \cap E_\pm$ by selecting $e_{13}^\pm, e_{14}^\pm$
as positive respectively negative basis  for $T_{e_{12}^\pm} S^2$. 
One checks that this definition is independent of the choice of a positive
orthonormal basis $e_1, e_2, e_3, e_4$ for $\rel^4$.
Now we define $\grass: G_{4,2} \rightarrow S^2_+ \times S^2_{-}$ by
\begin{eqnarray}
\nonumber
\grass \left( \left[\frac{v + i w}{\sqrt{2}}\right] \right) & = &
\sqrt{2}\, \Big(\prog_+(v \wedge w) , \prog_-(v \wedge w)\Big)\\ 
\label{conf.grass.def}
& = & \frac{1}{\sqrt{2}}\, \Big( v \wedge w + *(v \wedge w) ,
	v \wedge w - *(v \wedge w) \Big), 
\end{eqnarray}
and put $\grass_\pm = \prog_\pm \circ \grass: G_{4,2} \rightarrow S^2_\pm$.
Clearly $\grass$ is well-defined, smooth and injective.

\begin{proposition} \label{conf.grass}
With respect to the Fubini-Study metric on $G_{4,2}$ and the product
metric on $S^2_+ \times S^2_{-}$, the map 
$\grass: G_{4,2} \rightarrow S^2_+ \times S^2_{-}$ defined by 
(\ref{conf.grass.def}) is diffeomorphic and isometric up to a 
factor, more precisely
\ee \label{conf.grass.metric}
\grass^* g_{S^2_+ \times S^2_{-}} = 4 g_{FS}.
\ef
Moreover, the K\"ahler form $\omega$ as defined in
\bcite{muell.sver} has on $G_{4,2}$ the representation
\ee \label{conf.grass.kaehler}
	\omega
	= (\grass_+^* vol_{S^2_+} + \grass_-^* vol_{S^2_{-}})/2,
\ef
where the sphere factors $S^2_\pm$ are oriented as above.
\end{proposition}
{\pr Proof:} 
For an orthonormal system $v,w \in \rel^4$,
we put $z = (v + i w) / \sqrt{2} \in S^7 \subseteq \com^4$
and check that $z, \bar z$ is a complex orthonormal system in $\com^4$.
Extending  $v,w$ to an orthonormal basis $v,w,\tau_1,\tau_2$ of $\rel^4$,
we note that $z, \bar z, \tau_1, \tau_2 \in \com^4$ is actually
a hermitian basis of $\com^4$. Now for $\alpha \in \com$ and $j=1,2$ we have
\dd
\frac{d}{dt} \sum_{k=0}^3
(z + t \alpha \tau_j)_k^2|_{t = 0}
= 2 \alpha \sum_{k=0}^3  z_k (\tau_j)_k = 0,
\df
since $\bar z$ is perpendicular to $\tau_j$ in $\com^4$. Thus 
if $\pi:S^7 \rightarrow \com P^3$, $\pi(z) = [z]$, denotes the 
Hopf projection, then by (\ref{conf.grass-def}) we see that
\dd
T_{\pro(z)} G_{4,2} = 
{\rm span}_\com \{D\pi(z)\tau_1, D\pi(z)\tau_2\}.
\df
Now $D\pi(z)iz = 0$, and by definition of the Fubini-Study metric 
the restriction of $D\pi(z)$ to the horizontal space ${\{z\}_\com^\perp}$
is an isometry onto  $T_{\pro(z)}(\com P^3)$.
In particular the four vectors $D\pi(z)\tau_j, D\pi(z)i\tau_j$ 
for $j = 1,2$ are an orthonormal basis of $T_{\pro(z)} G_{4,2}$.
We calculate
\begin{eqnarray*}
D \grass(\pi(z))\,D \pi(z) \tau_j & = &
\sqrt{2}\,\frac{d}{d \theta}\, (\grass \circ \pro)
	\frac{1}{\sqrt{2}} \Big( (\cos \theta) v + i w + (\sin \theta) \tau_j \Big)\\
&= & 2\, \frac{d}{d \theta}\,
	\big( (\cos \theta) v + (\sin \theta) \tau_j\big)  \wedge w\\
& = & ( \tau_j \wedge w + *(\tau_j \wedge w),
	\tau_j \wedge w - *(\tau_j \wedge w) ),
\end{eqnarray*}
and
\begin{eqnarray*}
D \grass(\pi(z))\,D \pi(z) i\tau_j & = &
\sqrt{2}\,\frac{d}{d \theta}\, (\grass \circ \pro)
\frac{1}{\sqrt{2}} \Big(v+i(\cos \theta) w + (\sin \theta) \tau_j \Big)\\
&= & 2\, \frac{d}{d \theta}\,
v \wedge \big( (\cos \theta) w + (\sin \theta) \tau_j\big)\\
& = & (v \wedge \tau_j + *(v \wedge \tau_j),
v \wedge \tau_j - *(v \wedge \tau_j) ).
\end{eqnarray*}
Writing $(v,w,\tau_1,\tau_2) =: (e_1,e_2,e_3,e_4)$ we see that
$D (\grass \circ \pro)(z)$ maps as follows:
\ee \label{conf.grass.diff}
\ad
	e_3 \mapsto \sqrt{2} (-e_{14}^+ , e_{14}^-), \quad
	e_4 \mapsto \sqrt{2} (e_{13}^+ , -e_{13}^-), \\
	i e_3 \mapsto \sqrt{2} (e_{13}^+ , e_{13}^-), \quad
	i e_4 \mapsto \sqrt{2} (e_{14}^+ , e_{14}^-).
\af
\ef
In particular, $D\grass(\pi(z))$ 
maps an orthonormal basis of $T_{\pro(z)} G_{4,2}$
to twice an orthonormal basis of $T_{\grass(\pro(z))}(S^2 \times S^2)$,
which proves (\ref{conf.grass.metric}). Furthermore, $\grass$ is a 
local diffeomorphism by the inverse function theorem, 
hence $\grass(G_{4,2}) \subseteq S^2 \times S^2\ $ is open.
As $\grass(G_{4,2})$ is compact, non-empty and $S^2 \times S^2$ is connected,
we obtain that $\grass$ is surjective. As we already saw that 
$\grass$ is injective, it is a global diffeomorphism.\\
\\
The K\"ahler form $\omega$ on $\com P^3$ is defined in \bcite{muell.sver} by
\dd
\omega(D\pi \cdot \xi, D \pi \cdot \eta)
	= 2 g_{FS}(D \pro \cdot \xi , D \pro \cdot i \eta)
	\quad \mbox{for } \xi, \eta, \in \{z\}_\com^\perp.
\df
In  $T_{e_{12}^\pm} S_{\pm}^2$, the rotation by $+\pi/2$ is given by 
$J_\pm e_{13}^\pm = \pm e_{14}^\pm$, whence
\dd
	vol_{S^2_\pm}(\xi,\eta)
	= g_{S^2_\pm}(\xi,J_\pm \eta)
	\quad \mbox{for } \xi, \eta \in T_{e_{12}^\pm} S^2_\pm.
\df
Using (\ref{conf.grass.diff}) we see that
$
D (\grass_\pm \circ \pro) \cdot  i\xi
= J_\pm D (\grass_\pm \circ \pro) \xi
$
for any $\xi \in span_\com\{ \tau_1 , \tau_2 \}$.
Together with (\ref{conf.grass.metric}),
we obtain for all $\xi, \eta \in span_\com\{ \tau_1 , \tau_2 \}$
\begin{eqnarray*}
(\grass_+^* vol_{S^2_+} + \grass_-^* vol_{S^2_{-}})(D \pro \cdot \xi,D \pro \cdot \eta)
&  = &  \sum_{\pm} vol_{S^2_\pm}\big(
	D (\grass_\pm \circ \pro) \cdot \xi,
	D (\grass_\pm \circ \pro) \cdot \eta \big)\\
& = & \sum_{\pm} g_{S^2_\pm}\big(
	D (\grass_\pm \circ \pro) \cdot \xi,
	J_\pm D (\grass_\pm \circ \pro) \cdot \eta\big)\\
& = & \sum_{\pm} g_{S^2_\pm}\big(
	D (\grass_\pm \circ \pro) \cdot \xi,
	D (\grass_\pm \circ \pro) \cdot i \eta \big)\\
& = & g_{S^2_+ \times S^2_{-}}\big(
	D \grass \cdot ( D \pro \cdot \xi),
	D \grass \cdot (D \pro \cdot i \eta) \big)\\
& = & (\grass^* g_{S^2_+ \times S^2_{-}}) (D\pi \cdot \xi, D\pi\cdot i \eta)\\
& = & 4 g_{FS}(D\pi \cdot \xi, D\pi \cdot i \eta)\\
& = &  2 \omega(D\pi \cdot \xi, D\pi \cdot \eta),
\end{eqnarray*}
and (\ref{conf.grass.kaehler}) follows.
\proof
Next for any immersion $f:\rel^2 \to \rel^n$ we introduce a modified Gau{\ss} map by
\ee \label{conf.grass.gauss}
	\varphi := \grass \circ G: \rel^2 \rightarrow S^2_+ \times S^2_{-},
\ef
and denote by 
$\varphi_\pm := \prog_\pm \circ \varphi:\rel^2 \rightarrow S^2_\pm$
its corresponding projections.

\begin{proposition} \label{conf.pro}
The pullback of the volume form on $S^2_\pm$ via $\varphi_\pm$ is given by
\begin{equation} \label{conf.pro.vol}
\varphi_\pm^* vol_{S^2_\pm} =  (K \pm  R)\, vol_g,\quad \mbox{ where }\quad
R = 2 \langle A^\circ_{11} \wedge A^\circ_{12},\nu_1 \wedge \nu_2 \rangle. 
\end{equation}
Here we use an oriented orthonormal basis $e_1,e_2$ on $\sur$, and an oriented orthonormal
basis $\nu_1,\nu_2$ of normal vectors along $f$. In particular we have 
\ee \label{conf.pro.jac}
|R| \leq \frac{1}{2}|A^\circ|^2 \quad \mbox{ and }\quad
J \varphi_\pm = |K \pm R| \leq |K| + \frac{1}{2} |A^\circ|^2.
\ef
\end{proposition}
{\pr Proof:} We may assume that $f$ is (locally) the inclusion map, 
writing $e_{1,2}$ instead of $Df \cdot e_{1,2}$; also we write $e_{3,4}$ 
for $\nu_{1,2}$. It is easy to check that the definition 
of $R$ is independent of the choice of the (oriented)
bases. We have $G = \pi((e_1+ie_2)/\sqrt{2})$, whence by (\ref{conf.grass.def})
$$
\varphi_\pm = \sqrt{2} \Pi_\pm (e_1 \wedge e_2),
$$
Differentiating and using 
$\langle D_{e_1}e_1,e_1 \rangle  = \langle D_{e_2}e_2,e_2 \rangle = 0$
we obtain
$$
D\varphi_\pm \cdot e_k = 
\sqrt{2}\, \Pi_\pm  \Big(A(e_1,e_k) \wedge e_2 + e_1 \wedge A(e_2,e_k)\Big).
$$
Writing $A_{ij} = \alpha_{ij} e_3 + \beta_{ij} e_4$ and expanding yields
\begin{eqnarray*}
D\varphi_\pm \cdot e_k & = & 
\sqrt{2}\,\Pi_{\pm}\Big(
\alpha_{1k} e_3 \wedge e_2 + \beta_{1k} e_4 \wedge e_2 
+ \alpha_{2k} e_1 \wedge e_3 + \beta_{2k} e_1 \wedge e_4 \Big)\\
& = & (\alpha_{2k} \pm \beta_{1k}) e_{13}^\pm + (\mp \alpha_{1k} + \beta_{2k}) e_{14}^\pm.
\end{eqnarray*}
Now $e_{13}^\pm, e_{14}^\pm$ is a positive 
respectively negative orthonormal basis for $T_{e_{12}^\pm} S^2$,
therefore
\dd
\det(D \varphi_\pm)
= \pm \det \left(
	\begin{array}{cc}
	\alpha_{21} \pm \beta_{11}
	& \alpha_{22} \pm \beta_{12} \\
	\mp \alpha_{11} + \beta_{21}
	& \mp \alpha_{12} + \beta_{22} \\
	\end{array}
	\right).
\df
By choice of the bases at a point we can assume that
$\alpha_{12} = 0$ and $\mean = H e_3$ for 
$H = \alpha_{11} + \alpha_{22}$. Then $\beta_{11} = - \beta_{22}=: \beta$,
and $K = (H^2 - |A|^2)/2 = \alpha_{11} \alpha_{22} - \beta^2 - \beta_{12}^2$.
Hence
\begin{eqnarray*}
\det(D \varphi_\pm) & = & 
\det  \left(
\begin{array}{cc}
\beta & \pm \alpha_{22} + \beta_{12} \\
\mp \alpha_{11} + \beta_{12}
& - \beta 
\end{array}
\right)\\
& = & \alpha_{11}\alpha_{22} - \beta^2  - \beta_{12}^2 
\pm  (\alpha_{11}-\alpha_{22}) \beta_{12} \\
& = & K \pm  (\alpha_{11}-\alpha_{22}) \beta_{12}. 
\end{eqnarray*}
On the other hand from $A^\circ_{ij} = A_{ij} - \frac{1}{2} \mean g_{ij}$ we see that
$$
A^\circ_{11} \wedge A^\circ_{12} = 
\Big((\alpha_{11} - \frac{1}{2}H)\,e_3 + \beta_{11}e_4\Big) \wedge \beta_{12}e_4 
= \frac{1}{2} (\alpha_{11}-\alpha_{22})  \beta_{12}\, e_3 \wedge e_4.
$$
This proves (\ref{conf.pro.vol}), and (\ref{conf.pro.jac}) follows easily.
\proof

{\flushleft {\pr Proof of Theorem \ref{conf.flat} for n = 4: }} 
We have $(\varphi_{+}^\ast vol_{S^2_+} + \varphi_{-}^\ast vol_{S^2_{-}})/2 = K vol_g$
from (\ref{conf.grass.kaehler}) and (\ref{conf.grass.gauss}), as well as 
$|D\varphi|^2 = 4 |DG|^2 = 2 |A|^2$ by (\ref{conf.grass.metric}) and 
\bcite{muell.sver} 2.3. Recalling the discussion for $n=3$, it is therefore 
sufficient to find a solution $v:\rel^2 \to \rel$ of 
\ee \label{conf.flat4.equ}
- \Delta v = *(\varphi_+^* vol_{S^2_+} + \varphi_-^* vol_{S^2_{-}})/2\,
\mbox{ on }\, \rel^2 \quad \mbox{ and } \quad \lim_{z \rightarrow \infty} v(z) = 0,
\ef
which satisfies the estimates
\ee \label{conf.flat4.esti}
\|v\|_{L^\infty(\rel^2)},\,
\|D v\|_{L^2(\rel^2)},\,
\|D^2 v\|_{L^1(\rel^2)}
\leq C(\delta) \int_{\rel^2}|D \varphi|^2 \d \Lt.
\ef
Using (\ref{conf.pro.vol}), (\ref{conf.flat.curv}), (\ref{conf.pro.jac}) and
(\ref{conf.flat.curv4}), we obtain the following estimates, assuming without
loss of generality that both inequalities in (\ref{conf.flat.curv4}) are strict,
\ee \label{conf.flat4.ass}
\int_{\rel^2} J \varphi_\pm \d \mu_g  < 8 \pi - \delta \quad  \mbox{ and}
\quad \left| \int_{\rel^2}  \varphi_\pm^* vol_{S^2_\pm} \right| < 4 \pi.
\ef
As explained in \bcite{muell.sver} Proposition 3.4.1, we may assume 
using approximation that $\varphi$ is smooth and constant outside 
a compact set, while keeping the assumptions (\ref{conf.flat4.ass}).
Here, we do not assume anymore that $\varphi$ is obtained as the 
Gau{\ss} map of some surface. Our argument will essentially follow
\bcite{muell.sver} 3.4 and 3.5.\\
\\
Considering $\varphi_\pm$ as maps from $S^2$ to $S^2_\pm$ using the stereographic projection, 
we compute
\dd
\int_{\rel^2} \varphi_\pm^* vol_{S^2_\pm} = 4 \pi\, {\rm deg}(\varphi_\pm) \in 4 \pi \ganz,
\df
hence we conclude from (\ref{conf.flat4.ass}) that
\ee \label{conf.flat4.deg}
	{\rm deg}(\varphi_\pm) = 0.
\ef
Defining $G = \grass^{-1} \circ \varphi:\rel^2 \to G_{4,2} \subseteq \com P^3$, we get
\dd
\int_{\rel^2} G^* \omega 
= \int_{\rel^2} (\varphi_+^* vol_{S^2_+} + \varphi_-^* vol_{S^2_{-}})/2
= 2 \pi \Big( {\rm deg}(\varphi_+) + {\rm deg}(\varphi_-) \Big) = 0.
\df
Let $\pi:S^7 \to \com P^3$ be the Hopf projection. By Proposition 3.4.3  in \bcite{muell.sver}
the map $G$ has a lift $F:\rel^2 \to S^7$, i.e. $G = \pi \circ F$, whose Dirichlet integral 
is computed as follows, using $|DG|^2 = |D\varphi|^2/4$ and 
$G^\ast \omega = (\varphi_{+}^\ast vol_{S^2_+} + \varphi_{-}^\ast vol_{S^2_{-}})/2$,
\ee \label{conf.flat4.lift}
4 \int_{\rel^2} |D F|^2 \d \Lt
= \int_{\rel^2} |D \varphi|^2 \d \Lt
        + \| *(\varphi_+^* vol_{S^2_+} + \varphi_-^* vol_{S^2_{-}})\|_{W^{-1,2}(\rel^2)}^2.
\ef
Here for $w \in L^1_{loc}(\rel^2)$ the norm on the right hand side is 
\dd
\|w\|_{W^{-1,2}(\rel^2)} = \sup \Big\{ \int_{\rel^2} w \psi \d \Lt:
\psi \in C^\infty_0(\rel^2),\, \int_{\rel^2} |D \psi|^2 \d \Lt \leq 1\Big\}.
\df
By (\ref{conf.flat4.deg}), the number of preimages 
${\rm card\,}(\varphi_\pm^{-1}\{p\})$ must be even for almost every $p \in S^2_\pm$, 
whence (\ref{conf.flat4.ass}) implies
\dd
	vol_{S^2_\pm}(\varphi_\pm(\rel^2))
	\leq \frac{1}{2} \int_{\rel^2} J \varphi_\pm \d \mu_g
	< 4 \pi - \delta/2. 
\df
Therefore, we may choose open sets $U_\pm \subseteq S^2_\pm$ with 
$U_\pm \supseteq \varphi_\pm(\rel^2)$ and
\ee \label{conf.flat4.meas}
	vol_{S^2_\pm}(S^2_\pm - U_\pm) \geq \delta/2,
\ef
so that $\varphi(\rel^2) \subseteq U_+ \times U_- 
\subseteq S^2_+ \times S^2_-$. We shall now construct 
one-forms $\xi_\pm$ on $U_\pm$ with the properties
\begin{equation}
\label{conf.flat4.forms}
d\xi_\pm = vol_{S^2_\pm}|_{U_\pm} \quad \mbox{ and } \quad 
|\xi_\pm| \leq \frac{C}{\delta} \mbox{ on } U_\pm.
\end{equation}
Using euclidean coordinates $q = (x,y,z)$, we first define a one-form 
$\xi_{e_3}$ on $S^2 - \{e_3\}$ by 
$$
\xi_{e_3} = - \frac{x dy - y dx}{1-z} \quad \mbox{ or } \quad
\xi_{e_3}(q)\cdot v = - \frac{\langle e_3 \times q,v \rangle}{1-\langle e_3,q \rangle},
$$
where $\times$ denotes the cross product. In polar coordinates 
$x = \sin \vartheta \cos \varphi$, $y = \sin \vartheta \sin \varphi$ and 
$z = \cos \vartheta$, one readily checks that 
$\xi_{e_3} = - (1 + \cos \vartheta) d\varphi$
and hence 
$$
d\xi_{e_3} = \sin \vartheta\, d\vartheta \wedge \d\varphi = vol_{S^2} \quad
\mbox{ on } S^2 - \{e_3\},
$$ 
where $S^2$ is oriented by its exterior normal. Next for any $p \in S^2$ we
choose $T \in {\mathbb SO}(3)$ with $T p  = e_3$ and put 
$\xi_p = T^\ast \xi_{e_3}$ on $S^2 - \{p\}$. We have explicitely 
$$
\xi_p(q) \cdot v = - \frac{\langle p \times q,v \rangle }{1-\langle p,q \rangle}
\quad \mbox{on }S^2- \{p\},
$$
in particular
$$
d\xi_p = vol_{S^2} \, \mbox{ on }S^2-\{p\} \quad  \mbox{ and } \quad
|\xi_p(q)| \leq \frac{2}{|p - q|}.
$$
For $E \subseteq S^2$ closed with $vol_{S^2}(E) \geq \delta/2$ we now define 
on $U = S^2 - E$ the one-form
$$
\xi_E(q) = \mint{E} \xi_p(q) \d vol_{S^2}(p),
$$
which satisfies 
$$
d\xi_E = vol_{S^2}|_U \quad \mbox{ and } \quad 
|\xi_E(q)| \leq \frac{2}{\delta}  \int_{S^2} \frac{2}{|p-q|} \d vol_{S^2}(p) \leq \frac{C}{\delta}.
$$
The forms $\xi_\pm$ as in (\ref{conf.flat4.forms}) are obtained by choosing
orientation preserving isometries $T_\pm:S^2_\pm \to S^2$, and putting 
$\xi_\pm = T_\pm^\ast \xi_{E_\pm}$ where $E_\pm = T_\pm(S^2_\pm-U_\pm)$. 
Now define on $U_+ \times U_{-}$ the one-form $\xi = \Pi_+^\ast \xi_{+} +  \Pi_{-}^\ast \xi_{-}$,
and compute
$$
d \varphi^* \xi  = 
(\varphi_+^* vol_{S^2_+} + \varphi_-^* vol_{S^2_{-}})|_{U_+ \times U_-}, \quad \mbox{ and } \quad 
|\varphi^* \xi| \leq \frac{C}{\delta}\,|D \varphi|.
$$
As $\varphi_\pm(\rel^2) \subseteq U_\pm$, we can estimate 
for any $\psi \in C^\infty_0(\rel^2)$ 
\begin{eqnarray*}
\Big| \int_{\rel^2} *(\varphi_+^* vol_{S^2_{+}} + \varphi_-^* vol_{S^2_{-}}) \psi\,d\Lt \Big|
& = & \Big| \int_{\rel^2}  d(\varphi^* \xi)\, \psi \Big|\\
& = & \Big|\int_{\rel^2} (\varphi_+^\ast \xi_{+} + \varphi_{-}^\ast \xi_{-}) \wedge d\psi \Big|\\
& \leq & \frac{C}{\delta} \|D\varphi\|_{L^2(\rel^2)} \|D\psi\|_{L^2(\rel^2)},
\end{eqnarray*}
hence we get by the definition of the $W^{-1,2}$ norm and by (\ref{conf.flat4.lift})
\ee \label{conf.flat4.lift-esti}
\int_{\rel^2} |D F|^2 \d \Lt  \leq C(\delta) \int_{\rel^2} |D \varphi|^2 \d \Lt.
\ef
Now (\ref{conf.grass.kaehler}) and \bcite{muell.sver} 2.2 imply that
\dd
	\pro^* \grass^*
	(\prog_+^* vol_{S^2_+} + \prog_-^* vol_{S^2_-})/2
	= \pro^* \omega = \sum_{k=0}^3 i dz_k \wedge d \bar z_k.
\df
From $\varphi_\pm = \Pi_\pm \circ \varphi$ and $\varphi = \grass \circ \pro \circ F$
we therefore have 
$$
(\varphi_+^* vol_{S^2_+} + \varphi_-^* vol_{S^2_{-}})/2 
= F^* \sum_{k=0}^3 i dz_k \wedge d \bar z_k
= 2 \sum_{k=0}^3 \det(D F_k)\, dx \wedge dy.
$$
As in \bcite{muell.sver} Proposition 3.3.1, we apply \bcite{coi.lio.mey.semm}
to obtain the Hardy space estimate, combining with (\ref{conf.flat4.lift-esti}),
\dd
\| * (\varphi_{+}^* vol_{S^2_+} + \varphi_{-}^* vol_{S^2_-})/2 \|_{\har(\rel^2)}
\leq C \int_{\rel^2} |D F|^2 \d \Lt
\leq C(\delta) \int_{\rel^2} |D \varphi|^2 \d \Lt.
\df
Now \bcite{muell.sver} Theorem 3.2.1 yields the existence of a function
$v: \rel^2 \to \rel$ satisfying (\ref{conf.flat4.equ}) and
(\ref{conf.flat4.esti}), thereby proving the theorem also for $n = 4$.
\proof




\vspace{1cm}
\begin{tabbing}
{\sc Ernst Kuwert } \hspace{4.0cm}  \={\sc Reiner Sch\"atzle}\\
{\sc Mathematisches Institut} \> {\sc Mathematisches Institut}\\
{\sc Universit\"at Freiburg} \> {\sc Universit\"at T\"ubingen}\\
{\sc Eckerstra{\ss}e 1, D-79104 Freiburg} \>{\sc Auf der Morgenstelle 10,  D-72076 T\"ubingen}\\
ernst.kuwert@math.uni-freiburg.de  \> schaetz@everest.mathematik.uni-tuebingen.de
\end{tabbing}



\begin{thebibliography}{99}

\bibitem[BK03]{bauer.kuwert}
        {Bauer, M., Kuwert, E.,} {(2003)}
	{Existence of minimizing Willmore surfaces of prescribed genus},
	{International Mathematics Research Notices}
	{\bf {10}},
	{pp. 553--576}.
\bibitem[CLMS93]{coi.lio.mey.semm}
	{Coifman, R., Lions, P.L., Meyer, Y., Semmes, S.,} {(1993)}
	{Compensated compactness and Hardy spaces},
	{Journal de Math\'ematiques pures et appliqu\'ees}
	{\bf {72}},
	{pp. 247-286}.
\bibitem[GT]{gil.tru}
	{Gilbarg, D., Trudinger, N.S.,} {(1998)}
	{Elliptic Partial Differential Equations of Second Order},
	{third edition, Springer Verlag},
	{Berlin - Heidelberg - New York - Tokyo}.
\bibitem[HO82]{hoff.oss}
	{Hoffmann, D., Osserman, R.,} {(1982)}
	{The Area of the Generalized Gaussian Image
	and the Stability of Minimal Surfaces
	in $ S^n \mbox{ and } \rel^n $},
	{Mathematische Annalen}
	{\bf {260}},
	{pp. 437-452}.
\bibitem[Kli58]{klinge}
        {Klingenberg, W.,} {(1958)}
	{Contributions to Riemannian geometry in the large},
	{Annals of Mathematics}
	{\bf {69}},
	{pp. 654-666}.
\bibitem[Kus89]{kusner}
        {Kusner, R.,} {(1989)}
	{Comparison surfaces for the Willmore problem},
	{Pacific Journal of Mathematics}
	{\bf {138}},
	{pp. 317-345}.
\bibitem[KS04]{kuw.schae.will3}
	{Kuwert, E., Sch\"atzle, R.,} {(2004)}
	{Removability of point singularities of Willmore surfaces},
	{Annals of Mathematics}
	{\bf {160}},
	{pp. 315-357}.
\bibitem[KS07]{kuw.schae.will4}
        {Kuwert, E., Sch\"atzle, R.,} 
	{Minimizers of the Willmore functional with precribed conformal type},
        {Preprint 2007}.
\bibitem[LY82]{li.yau}
	{Li, P., Yau, S.T.,} {(1982)}
	{A new conformal invariant and its applications to
	the Willmore conjecture and the first
	eigenvalue on compact surfaces},
	{Inventiones Mathematicae}
	{\bf {69}},
	{pp. 269-291}.
\bibitem[MS95]{muell.sver}
	{M\"uller, S., Sverak, V.,} {(1995)}
	{On surfaces of finite total curvature},
	{Journal of Differential Geometry}
	{\bf {42}},
	{pp. 229-258}.
\bibitem[Sch02]{schmidt1}
        {Schmidt, M.,}
	{A proof of the Willmore conjecture},
        {arXiv:math/0203224v2 (2002)}.
\bibitem[Sim]{sim}
	{Simon, L.,} {(1983)}
	{Lectures on Geometric Measure Theory},
	{Proceedings of the Centre for Mathematical Analysis
	Australian National University},
	{Volume 3}.
\bibitem[Sim93]{sim.will}
	{Simon, L.,} {(1993)}
	{Existence of surfaces minimizing the
	Willmore functional},
	{Communications in Analysis and Geometry}
	{\bf {1}},
	{pp. 281-326}.
\bibitem[Tro]{tromba}
        {Tromba, A.J.,} {(1992)}
	{Teichm\"uller theory in Riemannian geometry},
	{ETH Lectures in Mathematics},
	{Birkh\"auser, Basel-Boston-Berlin}. 
\end{thebibliography}
\end{document}